\documentclass[runningheads]{llncs}
\usepackage{graphicx}

\usepackage{amsmath}
\usepackage{amsfonts}
\usepackage{amssymb}
\usepackage{array}
\usepackage[title]{appendix}
\usepackage[french,english]{babel}
\usepackage{capt-of}
\usepackage{chngcntr}
\usepackage{color}
\usepackage{easybmat}
\usepackage{lineno}
\usepackage{listings}
\usepackage{multirow,hhline}
\usepackage{nicefrac}
\usepackage[dvipsnames]{xcolor}

\renewcommand{\lstlistingname}{Snippet}

\makeatletter
\newcommand{\biggg}{\bBigg@{4}}
\newcommand{\Biggg}{\bBigg@{5}}
\makeatother

\newcommand{\dsp}{\displaystyle}
\newcommand{\dspfrac}[2]{\dsp{\frac{#1}{#2}}}
\renewcommand{\v}{\vspace{5mm}}
\newcommand{\ie}{\textit{id est} }
\newcommand{\eg}{\textit{exempli gratia} }

\newcommand{\ite}{\item[$\bullet$]}
\newcommand{\abs}{\text{abs}}
\newcommand{\rel}{\text{rel}}
\newcommand{\dt}{\delta t}
\newcommand{\dtref}{\dt_{\text{ref}}}
\newcommand{\tinit}{t^{\text{init}}}
\newcommand{\tend}{t^{\text{end}}}
\newcommand{\nsys}{n_{sys}}

\newcommand{\nstk}[1]{n_{st,#1}}
\newcommand{\nink}[1]{n_{in,#1}}
\newcommand{\noutk}[1]{n_{out,#1}}
\newcommand{\nintot}{\nink{tot}}
\newcommand{\nouttot}{\noutk{tot}}

\newcommand{\rnstk}[1]{\mathbb{R}^{\nstk{k}}}
\newcommand{\rnink}[1]{\mathbb{R}^{\nink{k}}}
\newcommand{\rnoutk}[1]{\mathbb{R}^{\noutk{k}}}
\newcommand{\rnintot}{\mathbb{R}^{\nink{tot}}}
\newcommand{\rnouttot}{\mathbb{R}^{\noutk{tot}}}
\newcommand{\Insys}{[\![1, \nsys]\!]}

\newcommand{\Inink}[1]{[\![1, \nink{#1}]\!]}
\newcommand{\Inoutk}[1]{[\![1, \noutk{#1}]\!]}
\newcommand{\FofTime}[1]{L([\tinit, \tend], #1)}
\newcommand{\ib}{\bar{\imath}}
\newcommand{\jb}{\bar{\jmath}}
\newcommand{\ub}{\underline{\smash{u}}}
\newcommand{\yb}{\underline{\smash{y}}}

\AtBeginDocument{\counterwithout{lstlisting}{chapter}}
\begin{document}

\title{IFOSMONDI Co-simulation Algorithm with Jacobian-Free Methods in PETSc\thanks{Supported by organization Siemens Industry Software.}}

\titlerunning{IFOSMONDI-JFM Co-simulation Algorithm with PETSc}

\author{Yohan Eguillon\inst{1,2}\orcidID{0000-0002-9386-4646} \and
Bruno Lacabanne\inst{1}\orcidID{0000-0003-1790-3663} \and
Damien Tromeur-Dervout\inst{2}\orcidID{0000-0002-0118-8100}}

\authorrunning{Y. Eguillon et al.}

\institute{
Siemens Industry Software, Roanne, France
\email{\{yohan.eguillon,bruno.lacabanne\}@siemens.com}\\
\and
Institut Camille Jordan, Université de Lyon, UMR5208 CNRS-U.Lyon1, Villeurbanne, France\\
\email{\{yohan.eguillon,damien.tromeur-dervout\}@univ-lyon1.fr}\\
}

\maketitle


\begin{abstract}
IFOSMONDI iterative algorithm for implicit co-simulation of coupled physical systems (introduced by the authors in july 2019 during the Simultech conference, p.176-186) enables us to solve the non-linear coupling function while keeping the smoothness of interfaces without introducing a delay. Moreover, it automatically adapts the size of the steps between data exchanges among the systems according to the difficulty of the solving of the coupling constraint. The latter was solved by a fixed-point algorithm in the original implementation whereas this paper introduces the JFM version (standing for \textit{Jacobian-Free Methods}). Most implementations of Newton-like methods require a jacobian matrix which can be difficult to compute in the co-simulation context, except in the case where the interfaces are represented by a Zero-Order-Hold (ZOH). As far as IFOSMONDI coupling algorithm uses Hermite interpolation for smoothness enhancement (up to Third-Order-Hold), we propose hereafter a new formulation of the non-linear coupling function including both the values and the time-derivatives of the coupling variables. This formulation is well designed for solving the coupling through jacobian-free Newton type methods. Consequently, successive function evaluations consist in multiple simulations of the systems on a co-simulation time-step using rollback. The orchestrator-workers structure of the algorithm enables us to combine the PETSc framework on the orchestrator side for the non-linear Newton-type solvers with the parallel integrations of the systems on the workers side thanks to MPI processes. Different non-linear methods will be compared to one another and to the original fixed-point implementation on a newly proposed 2-systems academic test-case (mass-spring-damper type) with direct feedthrough on both sides.

\end{abstract}

\newpage

\section{Introduction}

The use of co-simulation is increasing in the industry as it enables to connect and simulate systems with given interfaces (input and output variables) without disclosing the expertise inside. Hence, modellers can provide system architects with virtual systems as black boxes since the systems are able to interact through their interfaces. Among these interactions, the minimal requirement are quite simple: a system should at least be able to read the inputs given by the other systems, to simulate its physics inside (most of the time thanks to an embedded solver), and to provide outputs of the simulation to the other systems.

Besides its black box aspect protecting the know-how, co-simulation also enables physic-based decomposition (one system can represent the hydraulic part of a modular model, another the mechanical part, a third one the electrical part, and so on) and/or dynamics-based decomposition (some systems isolate the stiff state variables so that they do not constraint all the other states anymore during the simulation). In other words, the co-simulation opens many doors thanks to the modular aspect of the models handled.

The co-simulation field of research nowadays focuses on the numerical methods and algorithms that can be used to process simulations of such modular models. From the simplest implementations (non-iterative Jacobi) to very advanced algorithms \cite{Gomes2018survey,Sadjina2017,Schierz2012,Sicklinger2014,Benedikt2013}, co-simulation methods have been developped in different fields, showing that the underlying problems to be tackled are not straightforward. Some arising problems could clearly be identified since the moment it has become a center of interest for researchers, such as the delay between the given inputs and the retrieved outputs of a system (corresponding to the so-called "co-simulation step" or "macro-step"), the instabilities that might occur as a consequence of this delay \cite{Viel2014}, the discontinuities produced at each communication \cite{Busch2019}, the error estimation (and the use of it in order to adapt the macro step size) \cite{Schierz2012}, the techniques to solve the so-called "constraint function" corresponding to the interface of the systems \cite{Kubler2000,Schweizer2015}, and so on. Many of these problems have been addressed in papers either proposing an analysis, a method to solve them, or both.

In our previous paper \cite{Eguillon2019}, an iterative method that satisfies the interfaces consistency while avoiding discontinuities at each macro-step was proposed and compared to well-established methods (non-iterative Jacobi, zero-order hold iterative co-simulation \cite{Kubler2000}, and non-iterative algorithm enhancing variables smoothness \cite{Busch2019}). This algorithm was based on a fixed-point iterative method. Its evolution, presented in this paper, is based on iterative methods that normally require jacobian matrix computation, yet we use their jacobian-free version. The name of this method is IFOSMONDI-JFM, standing for \textit{Iterative and Flexible Order, SMOoth and Non-Delayed Interfaces, based on Jacobian-Free Methods}. The enhancements it brings to the classical IFOSMONDI method enable to solve cases that could not be solved by this previous version. The integration of an easily modulable jacobian-free method to solve the constraint function will be presented. In particular, the software integration was made possible thanks to the PETSc framework, a library that provides modulable numerical algorithms. The interfacing between PETSc and the co-simulation framework dealing with the systems, interfaces and polynomial representations will be detailed.

\section{Formalism and notations}
\label{section:formalism_and_notations}

\subsection{A word on JFM accronym}
\label{subsection:a_word_on_jfm_accronym}

In the whole paper, the abbreviation JFM will denote jacobian-free versions of iterative methods that are designed to bring a given function (so-called \textit{callback}) to zero and that normally require the computation of the jacobian matrix of the callback function. In particular, a fixed-point method does not meet these criteria: it is not a JFM, contrary to matrix-free versions of the Newton method, the Anderson method \cite{Anderson1965} or the non-linear GMRES method \cite{Oosterlee2000}.

\subsection{General notations}
\label{subsection:general_notations}

In this paper, we will focus on the explicit systems. In other words, we will consider that every system in the co-simulation is a dynamical system corresponding to an ODE (Ordinary Differential Equation). The time-domain of the ODEs considered will be written $[\tinit, \tend[$, and the variable $t$ will denote the time.

Let's consider $\nsys\in\mathbb{N}^*$ systems are involved: we will use the index $k\in\Insys$ to denote the $k$\up{th} system, and $\nstk{k}$, $\nink{k}$, and $\noutk{k}$ will respectively denote the number of state variables, the number of inputs, and the number of outputs of system $k$.

The time-dependant vectors of states, inputs and outputs of system $k$ will respectively be written $x_k\in\FofTime{\rnstk{k}}$, $u_k\in\FofTime{\rnink{k}}$, and $y_k\in\FofTime{\rnoutk{k}}$ where $L(A, B)$ denotes the set of functions of domain $A$ and co-domain $B$. Finally, we can write the ODE form of the system $k$:

\begin{equation}
\label{eq:ODE_sys_k}
\left\{
\begin{array}{lcl}
	\dot{x}_k(t) & = & f_k(t, x_k(t), u_k(t)) \\
	y_k(t) & = & g_k(t, x_k(t), u_k(t))
\end{array}
\right.
\end{equation}

\noindent Let $\nintot$ and $\nouttot$ respectively be the total amount of inputs $\sum_{k=1}^{\nsys}\nink{k}$ and the total amount of outputs $\sum_{k=1}^{\nsys}\noutk{k}$.

The total inputs and the total outputs vectors are simply concatenations of input and output vectors of every system:

\begin{equation}
\label{eq:total_io_vectors}
\begin{array}{lclcl}
	\ub(t) & = & (u_1(t)^T, \cdots, u_{\nsys}(t)^T)^T & \in & \FofTime{\rnintot} \\
	\yb(t) & = & (y_1(t)^T, \cdots, y_{\nsys}(t)^T)^T & \in & \FofTime{\rnouttot} \\
\end{array}
\end{equation}

\noindent Finally, a tilde symbol $\tilde{}$ will be added to a functional quantity to represent an element of its codomain. \eg, $\yb\in L([t^{[N]}, t^{[N+1]}[, \mathbb{R})$, so we can use $\tilde{\yb}$ to represent an element of $\rnouttot$.

\newpage

\subsection{Extractors and rearrangement}
\label{subsection:extractors_and_rearrangement}

In order to easily switch from global to local inputs, extractors are defined. For $k\in\Insys$, the extractor $E_k^u$ is the matrix defined by \eqref{eq:extractor_matrix}.

\begin{equation}
\label{eq:extractor_matrix}
\begin{array}{rcccccl}
	E^u_k =\bigg( & 0 & \Big| & \Big( I_{\nink{k}} \Big) & \Big| & 0 & \bigg) \\
	& \underbrace{\hspace{20mm}} & & \underbrace{\hspace{20mm}} & & \underbrace{\hspace{20mm}} & \\
	\multicolumn{2}{l}{\nink{k} \times \sum_{l=1}^{k-1}\nink{l}} & & \nink{k} \times \nink{k} & & \multicolumn{2}{r}{\nink{k} \times \sum_{l=k+1}^{\nsys}\nink{l}} \\
\end{array}
\end{equation}

\noindent where $\forall n\in\mathbb{N},\ I_n$ denotes the identity matrix of size $n$ by $n$.

The extractors enable to extract the inputs of a given system from the global inputs vector with a relation of the form $\tilde{u}_k = E^u_k \tilde{\ub}$. We have: $\forall k\in\Insys,\ E^u_k\in M_{\nink{k}, \nintot}(\{0, 1\})$.\\

A rearrangement operator will also be needed to handle concatenations of outputs and output derivatives. For this purpose, we will use the rearrangement matrix $R^{\yb}\in M_{\nouttot, \nouttot}(\{0, 1\})$ defined blockwise in \eqref{eq:rearrangement_matrix}.

\begin{equation}
\label{eq:rearrangement_matrix}
\begin{array}{c}
	R^{\yb} = \left(R^{\yb}_{K, L}\right)_{\substack{K\in[\![1,\ 2\ \nsys]\!]\\L\in[\![1,\ 2\ \nsys]\!]}} \\ \\
	\text{where} \\ \\
	R^{\yb}_{K, L} = \left\{
	\begin{array}{ll}
		I_{\noutk{K}} & \text{if}\ K \leqslant \nsys\ \text{and}\ L=2K-1 \\
		I_{\noutk{K-\nsys}} & \text{if}\ K > \nsys\ \text{and}\ L=2(K-\nsys) \\
		0 & \text{otherwise} \\
	\end{array}\right.
\end{array}
\end{equation}

\noindent The $R^{\yb}$ operator makes it possible to rearrange the outputs and output derivatives with a relation of the following form.

\begin{equation}
\label{eq:rearrangement_matrix_example}
\begin{array}{ccc}
	\left(\begin{array}{c}
	 	\\
		\tilde{\yb}\\ \\
		\hline\\
		\tilde{\dot{\yb}}\\ \\
	\end{array}\right)
	=
	\left(\begin{array}{c}
		\tilde{y}_1 \\
		\tilde{y}_2 \\
		\vdots \\
		\tilde{y}_{\nsys} \\
		\hline
		\tilde{\dot{y}}_1 \\
		\tilde{\dot{y}}_2 \\
		\vdots \\
		\tilde{\dot{y}}_{\nsys} \\
	\end{array}\right)
	=
	&
	\underbrace
	{
	\left(\begin{array}{ccccccc}
		\!I_{\noutk{1}}\!\! & 0 & 0 & 0 & \cdots & 0 & 0 \\
		0 & 0 & \!\!I_{\noutk{2}}\!\! & 0 & \cdots & 0 & 0 \\
		\vdots & \vdots & \vdots & \vdots & \ddots & \vdots & \vdots \\
		0 & 0 & 0 & 0 & \cdots & I_{\noutk{\nsys}}\!\!\! & 0 \\
		\hline
		0 & \!I_{\noutk{1}}\!\! & 0 & 0 & \cdots & 0 & 0 \\
		0 & 0 & 0 & \!\!I_{\noutk{2}}\! & \cdots & 0 & 0 \\
		\vdots & \vdots & \vdots & \vdots & \ddots & \vdots & \vdots \\
		0 & 0 & 0 & 0 & \cdots & 0 & \!\!\!I_{\noutk{\nsys}}\! \\
	\end{array}\right)
	}
	&
	\left(\begin{array}{c}
		\tilde{y}_1 \\
		\tilde{\dot{y}}_1 \\
		\tilde{y}_2 \\
		\tilde{\dot{y}}_2 \\
		\vdots \\
		\tilde{y}_{\nsys} \\
		\tilde{\dot{y}}_{\nsys} \\
	\end{array}\right)
	\\
	& R^{\yb} & \\
\end{array}
\end{equation}

\newpage

\subsection{Time discretization}
\label{subsection:time_discretization}

In the context of co-simulation, the $g_k$ and $f_k$ functions in \eqref{eq:ODE_sys_k} are usually not available directly. Thus, several co-simulation steps, so-called "macro-steps", are made between $\tinit$ and $\tend$. Let's introduce the notations of the discrete version of the quantities introduced in \ref{subsection:general_notations}.

A macro step will be defined by its starting and ending times, respectively denoted as $[t^{[N]}, t^{[N+1]}]$ for the $N$\up{th} macro-step. The macro-steps define a partition of the time-domain.

\begin{equation}
\label{eq:time_domain_partition}
\left\{
\begin{array}{lcl}
	[\tinit, \tend[ & = & \dsp{\bigcup_{N=0}^{N_{\max}-1}} [t^{[N]}, t^{[N+1]}[ \\ \\
	t^{[0]} & = & \tinit \\
	t^{[N_{\max}]} & = & \tend \\
	\multicolumn{3}{c}
	{
		\forall	N \in [\![0, N_{\max}-1]\!],\ t^{[N+1]} > t^{[N]}
	}
\end{array}
\right.
\end{equation}

\noindent Let $\delta t^{[N]}$ denote the size of the $N$\up{th} macro-step:

\begin{equation}
\label{eq:macro_step_size}
\left\{
\begin{array}{l}
	\forall N \in [\![0, N_{\max}-1]\!],\ \delta t^{[N]} = t^{[N+1]}-t^{[N]} > 0 \\ \\
	\dsp{\sum_{N=0}^{N_{\max}-1}} \delta t^{[N]} = \tend - \tinit
\end{array}
\right.
\end{equation}

\noindent Let $\mathbb{T}$ denote the set of possible macro-steps.

\begin{equation}
\label{eq:steps_set}
\mathbb{T} \overset{\Delta}{=} \{[a, b[\ |\ \tinit \leqslant a < b \leqslant \tend \} 
\end{equation}

\noindent On a given macro-step $[t^{[N]}, t^{[N+1]}[$, $N\in [0, N_{\max}]$, for all systems, the restrictions of the piecewise equivalents of $u_k$ and $y_k$ will be denoted by $u_k^{[N]}$ and $y_k^{[N]}$ respectively. In case several iterations are made on the same step, we will refer to the functions by a left superscript index $m$. Finally, we will denote the coordinate of these vectors with an extra subscript index.

\begin{equation}
\label{eq:discrete_io_funcs}
\renewcommand{\arraystretch}{1.5}
\begin{array}{llll}
	\multicolumn{4}{l}
	{
		\forall k\in\Insys,\ \forall N\in[\![0, N_{\max}]\!],\ \forall m\in[0, m_{\max}(N)],
	} \\
	\hspace{5mm} &
		\forall j\in\Inink{k},\ &
			^{[m]}u_{k, j}^{[N]} &
				\in L([t^{[N]}, t^{[N+1]}[, \mathbb{R}) \\
	\hspace{5mm} &
		\forall i\in\Inoutk{k},\ &
			^{[m]}y_{k, i}^{[N]} &
				\in L([t^{[N]}, t^{[N+1]}[, \mathbb{R}) \\
\end{array}
\renewcommand{\arraystretch}{1}
\end{equation}

\noindent In \eqref{eq:discrete_io_funcs}, $m_{\max}(N)$ denotes the number of iterations (minus one) done on the $N$\up{th} macro-step. $m_{max}(N)$ across $N$ can be plotted in order to see where the method needed to proceed more or less iterations.\\

All derived notations introduced in this subsection can also be applied to the total input and output vectors.

\begin{equation}
\label{eq:discrete_tot_io_funcs}
\renewcommand{\arraystretch}{1.5}
\begin{array}{llll}
	\multicolumn{4}{l}
	{
		\forall N\in[\![0, N_{\max}]\!],\ \forall m\in[\![0, m_{\max}(N)]\!],
	} \\
	\hspace{5mm} &
		\forall \jb\in\nintot,\ &
			^{[m]}\ub_{\jb}^{[N]} &
				\in L([t^{[N]}, t^{[N+1]}[, \mathbb{R}) \\
	\hspace{5mm} &
		\forall \ib\in\nouttot,\ &
			^{[m]}\yb_{\ib}^{[N]} &
				\in L([t^{[N]}, t^{[N+1]}[, \mathbb{R}) \\
\end{array}
\renewcommand{\arraystretch}{1}
\end{equation}

\noindent Indices $\ib$ and $\jb$ in \eqref{eq:discrete_tot_io_funcs} will be called \textit{global indices} in opposition to the \textit{local indices} $i$ and $j$ in \eqref{eq:discrete_io_funcs}.

\subsection{Step function}
\label{subsection:step_function}

Let $S_k,\ k\in\Insys$ be the \textit{ideal step function} of the $k$\up{th} system, that is to say the function which takes the system to its future state one macro-step forward.

\begin{equation}
\label{eq:step_function_ideal}
S_k:\left\{
\begin{array}{lcl}
	\mathbb{T} \times \FofTime{\rnink{k}} \times \rnstk{k} & \rightarrow & \rnoutk{k} \times \rnstk{k} \\
	(\tau,\ u_k,\ \tilde{x}) & \mapsto & S_k(\tau,\ u_k,\ \tilde{x})
\end{array}
\right.
\end{equation}

\noindent In practice, the state vector $\tilde{x}$ will not be explicited. Indeed, it will be embedded inside of the system $k$ and successive calls will either be done:

\begin{itemize}
\ite with $\tau$ beginning where the $\tau$ at the previous call of $S_k$ ended (moving on),
\ite with $\tau$ beginning where the $\tau$ at the previous call of $S_k$ started (step replay),
\ite with $\tau$ of the shape $[\tinit, t[$ with $t\in]\tinit, \tend[$ (first step).
\end{itemize}

Moreover, the $u_k$ argument only needs to be defined on the domain $\tau$ (not necessary on $[\tinit, \tend[$). Thus, $S_k$ will not be considered in the method, but the $\hat{S}_k$ function (\textit{practical step function)} defined hereafter will be considered instead. Despite $\hat{S}_k$ is not properly mathematically defined (the domain depends on the value of one of the arguments: $\tau$ and some quantities are hidden: the states), it does not lead to any problem, considering the hypotheses above.

\begin{equation}
\label{eq:step_function_pratical}
\begin{array}{c}
	\hat{S}_k:\left\{
	\begin{array}{lcl}
		\mathbb{T} \times L(\tau, \rnink{k}) & \mapsto & \rnoutk{k} \\
		(\tau, u_k) & \mapsto & \hat{S}_k(\tau, u_k)
	\end{array}
	\right.\\ \\
	\text{satisfying} \\ \\
	\hat{S}_k([t^{[N]}, t^{[N+1]}[,\ ^{[m]}u_k^{[N]}) =\ ^{[m]}y_k^{[N]}(t^{[N+1]})
\end{array}
\end{equation}

\noindent The $\hat{S}_k$ function is the one available in practice, namely in the FMI (Functionnal Mock-up Interface) standard.

\subsection{Extended step function}
\label{subsection:extended_step_function}

The values of the output variables might not be sufficient for every co-simulation scheme. It is namely the case for both classical IFOSMONDI and IFOSMONDI-JFM. Indeed, the time-derivatives of the outputs are also needed.

Let $\hat{\hat{S}}_k$ be the extension of $\hat{S}_k$ returning both the output values and derivatives.

\begin{equation}
\label{eq:step_function_pratical_extended}
\begin{array}{c}
	\hat{\hat{S}}_k:\left\{
	\begin{array}{lcl}
		\mathbb{T} \times L(\tau, \rnink{k}) & \mapsto & \rnoutk{k} \times \rnoutk{k} \\
		(\tau, u_k) & \mapsto & \hat{\hat{S}}_k(\tau, u_k)
	\end{array}
	\right.\\ \\
	\text{satisfying} \\ \\
	\hat{\hat{S}}_k([t^{[N]}, t^{[N+1]}[,\ ^{[m]}u_k^{[N]}) =
		\left(
			\begin{array}{c}
				^{[m]}y_k^{[N]}(t^{[N+1]}) \\ \\
				\dspfrac{d\ ^{[m]}y_k^{[N]}}{dt}(t^{[N+1]})
			\end{array}
		\right)
\end{array}
\end{equation}

\noindent If the derivatives are not available in practice, a finite difference approximation over $[t^{[N]}, t^{[N+1]}[$ can be made (see $\tilde{\hat{\hat{S}}}_k$ in \cite{Eguillon2019}).

\subsection{Connections}
\label{subsection:connections}

The connections between the systems will be denoted by a matrix filled with zeros and ones, with $\nouttot$ rows and $\nintot$ columns denoted by $\Phi$. Please note that if each output is connected to exactely one input, $\Phi$ is a square matrix. Moreover, it is a permutation matrix. Otherwise, if an output is connected to several inputs, more than one $1$ appears at the corresponding row of $\Phi$. In any case, there can neither be more nor less than one $1$ on each column of $\Phi$ considering that an input can neither be connected to none nor several outputs.

\begin{equation}
\label{eq:Phi}
\forall \ib\in\nouttot,\ \forall \jb\in\nintot,\ 
\Phi_{\ib, \jb} = \left\{
\begin{array}{ll}
	1 & \text{if output $\ib$ is connected to input $\jb$} \\
	0 & \text{otherwise}
\end{array}
\right.
\end{equation}

\noindent The \textit{dispatching} will denote the stage where the inputs are generated from their connected inputs, using the connections represented by $\Phi$.

\begin{equation}
\label{eq:dispatching}
\tilde{\ub} = \Phi^T \tilde{\yb}
\end{equation}

\noindent The \textit{coupling function} \eqref{eq:coupling_function} will denote the absolute difference between corresponding connected variables in a total input vector and a total output vector. In other words, it represents the absolute error beween a total input vector and the dispatching of a total output vector. The $\lambda$ subscript does not correspond to a quantity, it is a simple notation inherited from a Lagrange multipliers approach of systems coupling \cite{Schweizer2015}.

\begin{equation}
\label{eq:coupling_function}
g_{\lambda}:\left\{
\begin{array}{lcl}
	\rnintot \times \rnouttot & \rightarrow & \rnintot \\
	(\tilde{\ub}, \tilde{\yb}) & \mapsto & |\tilde{\ub} - \Phi^T \tilde{\yb}|
\end{array}
\right.
\end{equation}

\noindent The \textit{coupling condition} \eqref{eq:coupling_condition} is the situation where every output of the total output vector corresponds to its connected input in the total input vector.

\begin{equation}
\label{eq:coupling_condition}
g_{\lambda}(\tilde{\ub}, \tilde{\yb}) = 0_{\rnintot}
\end{equation}

\section{IFOSMONDI-JFM method}
\label{section:ifosmondi_jfm_method}

\subsection{Modified extended step function}
\label{subsection:modified_extended_step_function}

As in classical IFOSMONDI \cite{Eguillon2019}, the IFOSMONDI-JFM method preserves the $C^1$ smoothness of the interface variables at the communication times\\\noindent $(t^{[N]})_{N\in[\![1, N_{\max}-1]\!]}$. Thus, when a time $t^{[N]}$ has been reached, the input functions for every system will all satisfy the following property:

\begin{equation}
\label{eq:left_C1_smoothness}
\begin{array}{llcl}
	\multicolumn{4}{l}
	{
		\forall k\in\Insys,\ \forall m\in[\![0, m_{\max}(N)]\!],
	} \\ \\
	\hspace{5mm} \multirow{2}{*}{\biggg\{}
		& ^{[m]}u_k^{[N]}(t^{[N]}) & = & \ ^{[m_{\max}(N-1)]}u_k^{[N-1]}(t^{[N]}) \\ \\
		& \dspfrac{d\ ^{[m]}u_k^{[N]}}{dt}(t^{[N]}) & = & \ \dspfrac{d\ ^{[m_{\max}(N-1)]}u_k^{[N-1]}}{dt}(t^{[N]}) \\
\end{array}
\end{equation}

\noindent The IFOSMONDI-JFM also represents the inputs as $3$\up{rd} order polynomial (maximum) in order to satisfy the smoothness condition \eqref{eq:left_C1_smoothness} and to respect imposed values and derivatives at $t^{[N+1]}$ for every macro-step.

Knowing these constraints, it is possible to write a specification of the practical step function $\hat{\hat{S}}_k$ in the IFOSMONDI-JFM case (also applicable in the classical IFOSMONDI method):

\begin{equation}
\label{eq:step_function_specific}
\zeta_k:\left\{
\begin{array}{lcl}
	\mathbb{T} \times \rnink{k} \times \rnink{k} & \mapsto & \rnoutk{k} \times \rnoutk{k} \\
	(\tau, \tilde{u}_k, \tilde{\dot{u}}_k) & \mapsto & \zeta_k(\tau, \tilde{u}_k, \tilde{\dot{u}}_k)
\end{array}
\right.\\
\end{equation}

\noindent where the three cases discussed in \ref{subsection:step_function} have to be considered.

\subsubsection{Case $1$: Moving on:}
\label{subsubsection:moving_on}

In this case, the last call to $\zeta_k$ was done with a $\tau\in\mathbb{T}$ ending at current $t^{[N]}$. In other words, the system $k$ "reached" time $t^{[N]}$. The inputs were, at this last call: $^{[m_{\max}(N-1)]}u_k^{[N-1]}$.

To reproduce a behavior analog to the classical IFOSMONDI method, the inputs $^{[0]}u_k^{[N]}$ will be defined as the $2$\up{nd} order polynomial (or less) satisfying the three following constraints:

\begin{equation}
\label{eq:inputs_def_case1_moving_on}
\begin{array}{lcl}
	^{[0]}u_k^{[N]}(t^{[N]}) & = & ^{[m_{\max}(N-1)]}u_k^{[N-1]}(t^{[N]}) \\ \\
	\dspfrac{d\ ^{[0]}u_k^{[N]}}{dt}(t^{[N]}) & = & \dspfrac{d\ ^{[m_{\max}(N-1)]}u_k^{[N-1]}}{dt}(t^{[N]}) \\ \\
	^{[0]}u_k^{[N]}(t^{[N+1]}) & = & ^{[m_{\max}(N-1)]}u_k^{[N-1]}(t^{[N]}) \\
\end{array}
\end{equation}

\noindent The two first constraints guarantee the smoothness property \eqref{eq:left_C1_smoothness}, and the third one minimizes the risk of out-of-range values (as in classical IFOSMONDI method).

In this case, $\zeta_k$ in \eqref{eq:step_function_specific} is defined by the specification \eqref{eq:step_function_specific_props_case1}.

\begin{equation}
\label{eq:step_function_specific_props_case1}
\begin{array}{r}
	\zeta_k([t^{[N]}, t^{[N+1]}[, \cdot, \cdot) = \hat{\hat{S}}_k([t^{[N]}, t^{[N+1]}[,\ \underbrace{^{[0]}u_k^{[N]}}) \\
	\text{computed with \eqref{eq:inputs_def_case1_moving_on}}
\end{array}
\end{equation}

\noindent $2$\up{nd} and $3$\up{rd} arguments of $\zeta_k$ are unused.

\subsubsection{Case $2$: Step replay:}
\label{subsubsection:step_replay}

In this case, the last call to $\zeta_k$ was done with a $\tau\in\mathbb{T}$ starting at current $t^{[N]}$. In other words, the system did not manage to reach the ending time of the previous $\tau$ (either because the method did not converge, or because the step has been rejected, or another reason).

Two particular subcases have to be considered here: either the step we are computing is following the previous one in the iterative method detailed after this section, or the previous iteration has been rejected and we are trying to re-integrate the step starting from $\tau$ with a smaller size $\delta t^{[N]}$.

\paragraph{Subcase $2$.$1$: Following a previous classical step:}
\label{paragraph:following_a_previous_classical_step}

In this subcase, the last call of $\zeta_k$ was not only done with the same starting time, but also with the same step ending time $t^{[N+1]}$. The inputs were, at this last call: $^{[m-1]}u_k^{[N]}$ with $m\geqslant 1$, and satisfied the two conditions at $t^{[N]}$ of \eqref{eq:step_function_specific_props_case1}.

The jacobian-free iterative method will ask for given input values $\tilde{u}_k$ and time-derivatives $\tilde{\dot{u}}_k$ that will be used as constraints at $t^{[N+1]}$, thus $^{[m]}u_k^{[N]}$ will be defined as the $3$\up{rd} order polynomial (or less) satisfying the four following constraints:

\begin{equation}
\label{eq:inputs_def_case2_step_replay_subcase1}
\begin{array}{lclcl}
	^{[m]}u_k^{[N]}(t^{[N]})
		& = & ^{[m_{\max}(N-1)]}u_k^{[N-1]}(t^{[N]})
			& = & ^{[m-1]}u_k^{[N]}(t^{[N]}) \\ \\
	\dspfrac{d\ ^{[m]}u_k^{[N]}}{dt}(t^{[N]})
		& = & \dspfrac{d\ ^{[m_{\max}(N-1)]}u_k^{[N-1]}}{dt}(t^{[N]})
			& = & \dspfrac{d\ ^{[m-1]}u_k^{[N]}}{dt}(t^{[N]}) \\ \\
	^{[m]}u_k^{[N]}(t^{[N+1]}) & = & \tilde{u}_k \\ \\
	\dspfrac{d\ ^{[m]}u_k^{[N]}}{dt}(t^{[N+1]}) & = & \tilde{\dot{u}}_k \\
\end{array}
\end{equation}

\noindent The two firsts constraints ensure the \eqref{eq:left_C1_smoothness} smoothness property, and the third and fourth ones will enable the iterative method to find the best values and derivatives to satisfy the coupling condition.

In this subcase, $\zeta_k$ in \eqref{eq:step_function_specific} is defined by the specification \eqref{eq:step_function_specific_props_case2_subcase1}.

\begin{equation}
\label{eq:step_function_specific_props_case2_subcase1}
\begin{array}{r}
	\zeta_k([t^{[N]}, t^{[N+1]}[, \tilde{u}_k, \tilde{\dot{u}}_k) = \hat{\hat{S}}_k([t^{[N]}, t^{[N+1]}[,\ \underbrace{^{[m]}u_k^{[N]}}) \\
	\text{computed with \eqref{eq:inputs_def_case2_step_replay_subcase1}}
\end{array}
\end{equation}

\paragraph{Subcase $2$.$2$: Re-integrate a step starting from $t^{[N]}$ but with different $\delta t^{[N]}$ than at the previous call of $\zeta_k$:}
\label{paragraph:re_integrate_a_step_starting_from_t_N}

In this subcase, current $t^{[N+1]}$ is different from $\sup{(\tau)}$ with $\tau$ being the one used at the last call of $\zeta_k$.

As it shows that a step rejection just occured, we will simply do the same than in case $1$, as if we were moving on from $t^{[N]}$. In other words, all calls to $\zeta_k$ with $\tau$ starting at $t^{[N]}$ are "forgotten".

Please note that $^{[m_{\max}(N-1)]}u_k^{[N-1]}(t^{[N]})$ and $\dspfrac{d\ ^{[m_{\max}(N-1)]}u_k^{[N-1]}}{dt}(t^{[N]})$ can be retreived using the values and derivatives constraints at $t^{[N]}$ of the inputs at the last call of $\zeta_k$ thanks to the smoothness constraint \eqref{eq:left_C1_smoothness}.

\subsubsection{Case $3$: First step:}
\label{subsubsection:first_step}

In this particular case, we will do the same as in the other cases, except that we won't impose any constraint for the time-derivative at $\tinit$. That is to say:

\begin{itemize}
\ite at the first call of $\zeta_k$, we have $N=m=0$, we will only impose $^{[0]}u_k^{[0]}(\tinit) =\ ^{[0]}u_k^{[0]}(t^{[1]}) = u_{k}^{\text{init}}$ to have a zero order polynomial satisfying the initial conditions $u_{k}^{\text{init}}$ (supposed given),
\ite at the other calls, case $2$ will be used without considering the constraints for the derivatives at $\tinit$ (this will lower the polynomial's degrees). For \eqref{eq:inputs_def_case2_step_replay_subcase1}, the first condition becomes $^{[m]}u_k^{[N]}(\tinit) = u_{k}^{\text{init}}$, the second one vanishes, and the third ans fourth ones remain unchanged. For the subcase $2$.$2$, it can be considered that $^{[m_{\max}(-1)]}u_k^{[-1]}(\tinit) = u_{k}^{\text{init}}$, and $\dspfrac{d\ ^{[m_{\max}(-1)]}u_k^{[-1]}}{dt}(\tinit)$ will not be needed as it is a time-derivative in $\tinit$.
\end{itemize}

\v

Finally, we have $\zeta_k$ defined in every case, wrapping polynomial inputs computations and the integration done with $\hat{\hat{S}}_k$.

\subsection{Iterative method's callback function}
\label{subsection:iterative_method_s_callback_function}

The aim is to solve the co-simulation problem by using a jacobian-free version of an iterative method that usually requires a jacobian computation (see \ref{subsection:a_word_on_jfm_accronym}). Modern matrix-free versions of such algorithms make it possible to avoid perturbating the systems and re-integrating them for every input, as done in \cite{Schweizer2015}, in order to compute a finite-differences jacobian matrix. This saves a lot of integrations over each macro-step and gains time.

Nevertheless, on every considered macro-step $\tau$, a function to be brought to zero has to be defined. This so-called \textit{JFM's callback} (standing for \textit{Jacobian-Free Method's callback}) presented hereafter will be denoted by $\gamma_{\tau}$. In zero-order hold co-simulation, this function if often $\tilde{\ub}-\Phi^T\tilde{\yb}$ (or equivalent) where $\tilde{\yb}$ are the output at $t^{[N+1]}$ generated by constant inputs $\tilde{\ub}$ over $[t^{[N]}, t^{[N+1]}[$.

In IFOSMONDI-JFM, we will only enable to change the inputs at $t^{[N+1]}$, the smoothness condition at $t^{[N]}$ guaranteeing that the coupling condition \eqref{eq:coupling_condition} remains satisfied at $t^{[N]}$ if it was satisfied before moving on to the step $[t^{[N]}, t^{[N+1]}[$. The time-derivatives will also be considered in order to maintain $C^1$ smoothness, so the coupling condition \eqref{eq:coupling_condition} will also be applied to these time-derivatives.

Finally, the formulation of the JFM's callback for IFOSMONDI-JFM is:

\begin{equation}
\label{eq:jfm_s_callback}
\gamma_{\tau}:\left\{
\begin{array}{lcl}
	\rnintot \times \rnintot & \rightarrow & \rnintot \times \rnintot \\
	\left(\begin{array}{c}
		\tilde{\ub} \\
		\tilde{\dot{\ub}}
	\end{array}\right)
	& \mapsto &
	\left(\begin{array}{c}
		\tilde{\ub} \\
		\tilde{\dot{\ub}}
	\end{array}\right)
	-
	\left(\begin{array}{cc}
		\Phi^T & 0 \\
		0 & \Phi^T
	\end{array}\right)
	R^{\yb}
	\left(\begin{array}{c}
		\zeta_1\left(\tau, E^u_1 \tilde{\ub}, E^u_1 \tilde{\dot{\ub}}\right) \\
		\vdots \\
		\zeta_{\nsys}\left(\tau, E^u_{\nsys} \tilde{\ub}, E^u_{\nsys} \tilde{\dot{\ub}}\right) \\
	\end{array}\right)
\end{array}
\right.
\end{equation}

\subsubsection{Link with the fixed-point implementation:}

The formulation \eqref{eq:jfm_s_callback} can be used to represent the expression of the fixed-point $\Psi_{\tau}$ function. The latter has been introduced in classical IFOSMONDI algorithm \cite{Eguillon2019} where a fixed-point method was used instead of a JFM one. We can now rewrite a proper expression of $\Psi_{\tau}$ including the time-derivatives.

\begin{equation}
\label{eq:fixed_points_callback}
\Psi_{\tau}:\left\{
\begin{array}{lcl}
	\rnintot \times \rnintot & \rightarrow & \rnintot \times \rnintot \\
	\left(\begin{array}{c}
		\tilde{\ub} \\
		\tilde{\dot{\ub}}
	\end{array}\right)
	& \mapsto &
	\left(\begin{array}{c}
		\tilde{\ub} \\
		\tilde{\dot{\ub}}
	\end{array}\right)
	-
	\gamma_{\tau}
	(
	\left(\begin{array}{c}
		\tilde{\ub} \\
		\tilde{\dot{\ub}}
	\end{array}\right)
	)
	\\
	&&
	=
	\left(\begin{array}{cc}
		\Phi^T & 0 \\
		0 & \Phi^T
	\end{array}\right)
	R^{\yb}
	\left(\begin{array}{c}
		\zeta_1\left(\tau, E^u_1 \tilde{\ub}, E^u_1 \tilde{\dot{\ub}}\right) \\
		\vdots \\
		\zeta_{\nsys}\left(\tau, E^u_{\nsys} \tilde{\ub}, E^u_{\nsys} \tilde{\dot{\ub}}\right) \\
	\end{array}\right)
\end{array}
\right.
\end{equation}

\noindent $\Psi_{\tau}$ was refered as $\Psi$ in \cite{Eguillon2019} and did not include the derivatives in its formulation, yet the smoothness enhancement done by the Hermite interpolation led to an underlying use of these derivatives.

When the result of the $m$\up{th} iteration is available, a fixed-point iteration on macro-step $\tau=[t^{[N]}, t^{[N+1]}[$ is thus simply done by:

\begin{equation}
\label{eq:fixed_point_iteration}
\left(\begin{array}{c}
	^{[m+1]}\tilde{\ub} \\
	^{[m+1]}\tilde{\dot{\ub}}
\end{array}\right)
:=
\Psi_{\tau}
(
\left(\begin{array}{c}
	^{[m]}\tilde{\ub} \\
	^{[m]}\tilde{\dot{\ub}}
\end{array}\right)
)
\end{equation}

\subsection{First and last integrations of a step}
\label{subsection:first_and_last_iterations_of_a_step}

The first iteration of a given macro-step $\tau\in\mathbb{T}$ is a particular case to be taken into account. Considering the breakdown presented in subsection \ref{subsection:step_function}, this corresponds to case $1$, case $2$ subcase $2$.$2$, case $3$ first bullet point, and case $3$ second bullet point when falling into subcase $2$.$2$.

All these cases have something in common: they denote calls to $\zeta_k$ using a $\tau$ argument that has never been used in a previous call of $\zeta_k$. In these cases, the latter function is defined by \eqref{eq:step_function_specific_props_case1}.

For this reason, the first call of $\gamma_{\tau}$ for a given macro-step $\tau$ will be done before applying the JFM. Then, every time the JFM will call $\gamma_{\tau}$, the $(\zeta_k)_{k\in\Insys}$ functions called by $\gamma_{\tau}$ will behave the same way.

Once the JFM method ends, if it converged, a last call to $\gamma_{\tau}$ is made with the solution $\big((^{[m_{\max}(N)]}\tilde{\ub}^{[N]})^T,\ (^{[m_{\max}(N)]}\tilde{\dot{\ub}}^{[N]})^T\big)^T$ for the systems to be in a good state for the next step (as explained in subsection \ref{subsection:step_function}, the state of a system is hidden but affected at each call to a step function).

\subsection{Step size control}
\label{subsection:step_size_control}

The step size control is defined with the same \textit{rule-of-thumbs} than the one used in \cite{Eguillon2019}. The adaptation is not done on an error-based criterion such as in \cite{Schierz2012}, but with a predefined rule based on the convergence of the iterative method (yes/no).

A reference step size $\dtref \in \mathbb{R_*^+}$ is defined for any simulation with IFOSMONDI-JFM method. It will either act as initial macro-step size, and maximum step size. At some points, the method will be allowed to reduce this step in order to help the convergence of the JFM.

The convergence criterion for the iterative method is defined by the rule \eqref{eq:convergence_criterion}.

\begin{equation}
\label{eq:convergence_criterion}
\begin{array}{l}
	\text{Given}\ (\varepsilon_{\abs}, \varepsilon_{\rel})\in(\mathbb{R_+^*})^2, \\
	\text{convergence is reached when}\
		\left|
		\gamma_{\tau}
		\left(
		\begin{array}{c}
		\tilde{\ub} \\
		\tilde{\dot{\ub}}
		\end{array}
		\right)
		\right|
		<
		\left|
		\left(
		\begin{array}{c}
		\tilde{\ub} \\
		\tilde{\dot{\ub}}
		\end{array}
		\right)
		\right|
		\varepsilon_{\rel}
		+
		\left|
		\left(
		\substack{
		1 \\
		\vdots \\
		1
		}
		\right)
		\right|
		\varepsilon_{\abs}
\end{array}
\end{equation}

\noindent When the iterative method does not converge on the step $[t^{[N]}, t^{[N+1]}[$, either because a maximum number of iterations is reached or for any other reason (linear search does not converge, a Krylov internal method finds a singular matrix, ...), the step will be rejected and retried on the half \eqref{eq:step_size_def}. Otherwise, once the method converged on $[t^{[N]}, t^{[N+1]}[$, the next integration step $\tau$ tries to increase the size of $30\%$, without exceeding $\dtref$.

Once the iterative method exits on $\tau_{\text{old}}$, the next step $\tau_{\text{new}}$ is defined by:

\begin{equation}
\label{eq:step_size_def}
\begin{small}
\tau_{\text{new}} = \left\{
\begin{array}{lr}
	\multicolumn{2}{l}{\left[ \sup(\tau_{\text{old}}), \min\!\bigg\{\tend, \sup(\tau_{\text{old}}) + \max\!\Big\{\dtref, 1.3\ \big(\sup(\tau_{\text{old}})-\inf(\tau_{\text{old}})\big) \Big\}\bigg\}\right[
	} \\
	& \text{if convergence \eqref{eq:convergence_criterion} was reached} \\ \\
	\multicolumn{2}{l}{\left[\inf(\tau_{\text{old}}), \inf(\tau_{\text{old}}) + \dspfrac{\sup(\tau_{\text{old}})-\inf(\tau_{\text{old}})}{2}\right[} \\
	& \text{otherwise (divergence)} \\
\end{array}
\right.
\end{small}
\end{equation}

\noindent When $\varepsilon_{\abs} = \varepsilon_{\rel}$, these values will be denoted by $\varepsilon$.

\section{Note on the implementation}
\label{section:note_on_the_implementation}

Our implementation is based on an orchestrator-worker architecture, where $\nsys+1$ processes are involved. One of them is dedicated to global manipulations: the \textit{orchestrator}. It is not responsible of any system and only deals with global quantities (such as the time, the step $\tau$, the $\tilde{\ub}$ and $\tilde{\yb}$ vectors and the corresponding time-derivatives, and so on). The $\nsys$ remaining processes, the \textit{workers}, are reponsible of one system each. They only deal with local quantities related to the system they are responsible of.

\subsection{Parallel evaluation of $\gamma_{\tau}$ using MPI}
\label{subsection:parallel_evaluation_of_gamma_using_MPI}

An evaluation of $\gamma_{\tau}$ consists in evaluations of the $\nsys$ functions $(\zeta_k)_{k\in\Insys}$, plus some manipulations of vectors and matrices \eqref{eq:jfm_s_callback}. An evaluation of a single $\zeta_k$ for a given $k\in\Insys$ consists in polynomial computations and an integration \eqref{eq:step_function_specific_props_case1} \eqref{eq:step_function_specific_props_case2_subcase1} through a call of the corresponding $\hat{\hat{S}}_k$ function \eqref{eq:step_function_pratical_extended}.

A single call to $\gamma_{\tau}$ can be evaluated parallely by $\nsys$ processes, each of them carying out the integration of one of the systems. To achieve this, the MPI standard (standing for \textit{Message Passing Interface} has been used, as the latter provides routine to handle multi-process communications of data.

As the $k$\up{th} system only needs $E^u_k \tilde{\ub}$ and $E^u_k \tilde{\dot{\ub}}$ (see \eqref{eq:extractor_matrix}) among $\tilde{\ub}$ and $\tilde{\dot{\ub}}$, the data can be send in an optimized manner from an orchestrator process to $\nsys$ workers by using the \verb!MPI_Scatterv! routine.

Analogously, each worker process will have to communicate their contribution both to the outputs and their derivatives (assembling the block vector at the right of the expression \eqref{eq:jfm_s_callback}). This can be done by using the \verb!MPI_Gatherv! routine. 

Finally, global quantities such as $\tau$, $m$, the notifications of statuses and so on can be done easily thanks to the \verb!MPI_Broadcast! routine.

\subsection{Using PETSc for the JFM}
\label{subsection:using_petsc_for_the_jfm}

PETSc \cite{PetscWebPage,PetscParallelEfficiency} is a library used for parallel numerical computations. In our case, the several matrix-free versions of the Newton method and variants implemented in PETSc were very attractive. Indeed, the flexibility of this library at runtime enables the use of command-line arguments to control the resolution: \verb!-snes_mf! orders the use of a matrix-free non-linear solver, \verb!-snes_type newtonls!, \verb!anderson! \cite{Anderson1965} and \verb!ngmres! \cite{Oosterlee2000} are various usable solving methods that can be used as JFMs, \verb!-snes_atol!, \verb!-snes_rtol! and \verb!-snes_max_it! control the convergence criterion, \verb!-snes_converged_reason!, \\\verb!-snes_monitor! and \verb!-log_view! produce information and statistics about hte run, ...

This subsection proposes a solution to use these PETSc implementations in a manner that is compliant with the parallel evaluation of the JFM's callback \eqref{eq:jfm_s_callback}. This implementation has been used to generate the results of section \ref{section:results_on_a_test_cases}.

First of all, PETSc needs a view on the environment of the running code: the processes, and their relationships. In our case, the $\nsys+1$ processes of the orchestrator-worker architecture are not dedicated to the JFM. Thus, PETSc runs on the orchestrator process only. In terms of code, this can be done by creating PETSc objects referring to \verb!PETSC_COMM_SELF! communicator on the orchestrator process, and creating no PETSc object on the workers.

The callback $\gamma_{\tau}$ implements internally the communications with the workers, and is given to the PETSc \verb!SNES! object. The \verb!SNES! non-linear solver will call this callback blindly, and the workers will be triggered \textit{behind the scene} for integrations, preceded by the communications of the $\big((^{[m_{\max}(N)]}\tilde{\ub}^{[N]})^T,\ (^{[m_{\max}(N)]}\tilde{\dot{\ub}}^{[N]})^T\big)^T$ values asked by the \verb!SNES! and followed by the gathering of the outputs and related derivatives. The latters are finally returned to PETSc by the callback on the orchestrator side, after reordering and dispatching them as in \eqref{eq:jfm_s_callback}.

\subsection{JFM's callback implementation}
\label{subsection:jfm_s_callback_implementation}

In this section, a suggestion of implementation is proposed for the $\gamma_{\tau}$ function, both on the orchestrator side and on the workers side. Precisions about variables in the snippets are given below them.

The aim is not to show the code that has been used to generate the results of section \ref{section:results_on_a_test_cases}, but to figure out how to combine the PETSc and MPI standard (PETSc being based on MPI) to implement a parallel evaluation of $\gamma_{\tau}$.

By convention, the process of rank $0$ is the orchestrator, and any process of rank $k\in\Insys$ will be responsible of system $k$.\\

\vspace{2mm}

\begin{lstlisting}[caption=JFM's callback on the orchestrator side ($\gamma_{\tau}$), label=listing:gamma]
PetscErrorCode JFM_callback(SNES /* snes */, Vec u_and_du, Vec res, void *ctx_as_void)
{
   MyCtxType *ctx = (MyCtxType*)ctx_as_void;
   const int order = DO_A_STEP;
   PetscScalar const * pscalar_u_and_du;
   PetscScalar * pscalar_res;
   size_t k;

   // conversion PetscScalar -> C double
   VecGetArrayRead(u_and_du, &pscalar_u_and_du);
   for (k = 0; k < ctx->n_in_tot * 2; k++)
      ctx->double_u_and_du[k] = (double)(pswork_x[k]);
   VecRestoreArrayRead(u_and_du, &pscalar_u_and_du);

   // Notify workers that we want them to run,
   // and telling them what tau is
   MPI_Bcast(&order, 1, MPI_INT, 0, MPI_COMM_WORLD);
   MPI_Bcast(&(ctx->t_N),   1, MPI_DOUBLE, 0, MPI_COMM_WORLD);
   MPI_Bcast(&(ctx->t_Np1), 1, MPI_DOUBLE, 0, MPI_COMM_WORLD);

   // Apply extrators and communicate at the same time:
   // values, then derivatives
   MPI_Scatterv(ctx->double_u_and_du, ctx->in_sizes, ctx->in_offsets, MPI_DOUBLE, NULL, 0, MPI_DOUBLE, 0, MPI_COMM_WORLD);
   MPI_Scatterv(ctx->double_u_and_du + ctx->n_in_tot, ctx->in_sizes, ctx->in_offsets, MPI_DOUBLE, NULL, 0, MPI_DOUBLE, 0, MPI_COMM_WORLD);

   /* Workers proceed integration here */

   // Assemble vector R^{\bar{y}} * (\zeta_1^T, ... \zeta_2^T)^T directly
   // while communicating values and derivatives
   MPI_Gatherv_outputs(MPI_IN_PLACE, 0, MPI_DOUBLE, ctx->work1_n_out_tot, ctx->out_sizes, ctx->out_offsets, MPI_DOUBLE, 0, MPI_COMM_WORLD);
   MPI_Gatherv_outputs(MPI_IN_PLACE, 0, MPI_DOUBLE, ctx->work2_n_out_tot, ctx->out_sizes, ctx->out_offsets, MPI_DOUBLE, 0, MPI_COMM_WORLD);

   // Dispatching (equivalent of [[Phi^T, 0], [0, Phi^T]])
   dispatch(ctx->work1_n_out, ctx->out_sizes, ctx->n_sys,
      ctx->double_res, ctx->in_sizes, ctx->connections);
   dispatch(ctx->work2_n_out, ctx->out_sizes, ctx->n_sys,
      ctx->double_res + ctx->n_in_tot, ctx->in_sizes, ctx->connections);

   // Difference between original entries and permuted outputs
   for (k = 0; k < ctx->n_in_tot * 2; k++)
      ctx->double_res[k] = ctx->double_u_and_du[k] - ctx->double_res;

   // conversion C double -> PetscScalar
   VecGetArray(res, &pscalar_res);
   for (k = 0; k < ctx->n_in_tot * 2; k++)
      pswork_f[k] = (PetscScalar)(ctx->double_res[k]);
   VecRestoreArray(res, &pscalar_res);

   return 0;
}
\end{lstlisting}

\vspace{2mm}

\noindent In the code snippet \ref{listing:gamma}, the function \verb!JFM_callback! is the one that is given to the PETSc \verb!SNES! object with \verb!SNESSetFunction!. The context pointer \verb!ctx! can be anything that can be used to have access to extra data inside of this callback. The principle is: when \verb!SNESSolve! is called, the callback function which has been given to the \verb!SNES! object will be called an unknown number of times. For this example, we suggested a context structure \verb!MyCtxType! at least containing:

\begin{itemize}
\ite \verb!t_N!, \verb!t_Np1! the boundary times of $\tau$, \ie $t^{[N]}$ and $t^{[N+1]}$ (as \verb!double! each),
\ite \verb!n_in_tot! the total number of inputs $\nintot$ (as \verb!size_t!),
\ite \verb!double_u_and_du! an array dedicated to the storage of $(\tilde{\ub}^T, \tilde{\dot{\ub}}^T)^T$ (as \verb!double*!),
\ite \verb!in_sizes! the array containing the number of inputs for each process\\$(\nink{k})_{k\in[\![0, \nsys]\!]}$ including process $0$ (with the convention $\nink{0}=0$) (as \verb!int*!),
\ite \verb!in_offsets! the memory displacements $\left(\sum_{l=1}^{k}\nink{l}\right)_{k\in[\![0, \nsys]\!]}$ for inputs scattering for each process (as \verb!int*!),
\ite \verb!work1_n_out_tot! and \verb!work2_n_out_tot! two arrays of size $\nouttot$ for temporary storage (as \verb!double*!),
\ite \verb!out_sizes! and \verb!out_offsets! two arrays analogous to \verb!in_sizes! and \verb!in_offsets! respectively, considering the outputs,
\ite \verb!n_sys! tot number of systems $\nsys$ (as \verb!size_t!),
\ite \verb!double_res! an array of size $2\ \nintot$ dedicated to the storage of the result of $\gamma_{\tau}$ (as \verb!double*!), and
\ite \verb!connections! any structure to represent the connections between the systems $\Phi^T$ (a full matrix might be a bad idea as $\Phi$ is expected to be very sparse).
\end{itemize}

Finally, \verb!dispatch! is expected to be a function processing the dispatching \eqref{eq:dispatching} of the values given in its first argument into the array pointed by its fourth argument.\\

\newpage

On the workers side, the running code section is the one in the snippet \ref{listing:zeta}.

\begin{lstlisting}[caption=JFM's callback on the worker side ($\zeta_k$ and communications), label=listing:zeta]
/* ... */

while (1)
{
   MPI_Bcast(&order, 1, MPI_INT, 0, me__->comm);
   if (order != DO_A_STEP)
      break;

   // get tau
   MPI_Bcast(&t_N,   1, MPI_DOUBLE, 0, MPI_COMM_WORLD);
   MPI_Bcast(&t_Np1, 1, MPI_DOUBLE, 0, MPI_COMM_WORLD);

   // receive relevant inputs and derivatives for this system
   MPI_Scatterv(NULL, ctx->in_sizes, ctx->in_offsets, MPI_DOUBLE,
      sys_inputs, sys_n_in, MPI_DOUBLE, 0, me__->comm);
   MPI_Scatterv(NULL, ctx->in_sizes, ctx->in_offsets, MPI_DOUBLE,
      sys_dinputs, sys_n_in, MPI_DOUBLE, 0, me__->comm);

   /* integration: */
   zeta_do_a_step(t_N, t_Np1, inputs, sys_dinputs, // [in]
      sys_outputs, sys_doutputs); // [out]

   // send the outputs and derivatives (results of zeta)
   MPI_Gatherv_outputs(sys_outputs, sys_n_out, MPI_DOUBLE,
      NULL, NULL, NULL, MPI_DOUBLE, 0, MPI_COMM_WORLD);
   MPI_Gatherv_outputs(sys_doutputs, sys_n_out, MPI_DOUBLE,
      NULL, NULL, NULL, MPI_DOUBLE, 0, MPI_COMM_WORLD);
}

/* ... */
\end{lstlisting}

\noindent Please note that the orchestrator process has to explicitely send an order different from \verb!DO_A_STEP! (with \verb!MPI_Bcast!) to notify the workers that the callback will not be called anymore.

Nonetheless, this order might not be send right after the call to \verb!SNESSolve! on the orchestrator side. Indeed, if the procedure converged, a last call has to be made explicitely in the orchestrator (see \ref{subsection:first_and_last_iterations_of_a_step}).

An other explicit call to \verb!JFM_callback! should also be explicitly made on the orchestrator side before the call of \verb!SNESSolve! (as also explained in \ref{subsection:first_and_last_iterations_of_a_step}).\\

\newpage

Finally, figure \ref{fig:callback_sheme} presents a schematic view of these two snippets running parallely.

\begin{center}
\includegraphics[scale=0.25]{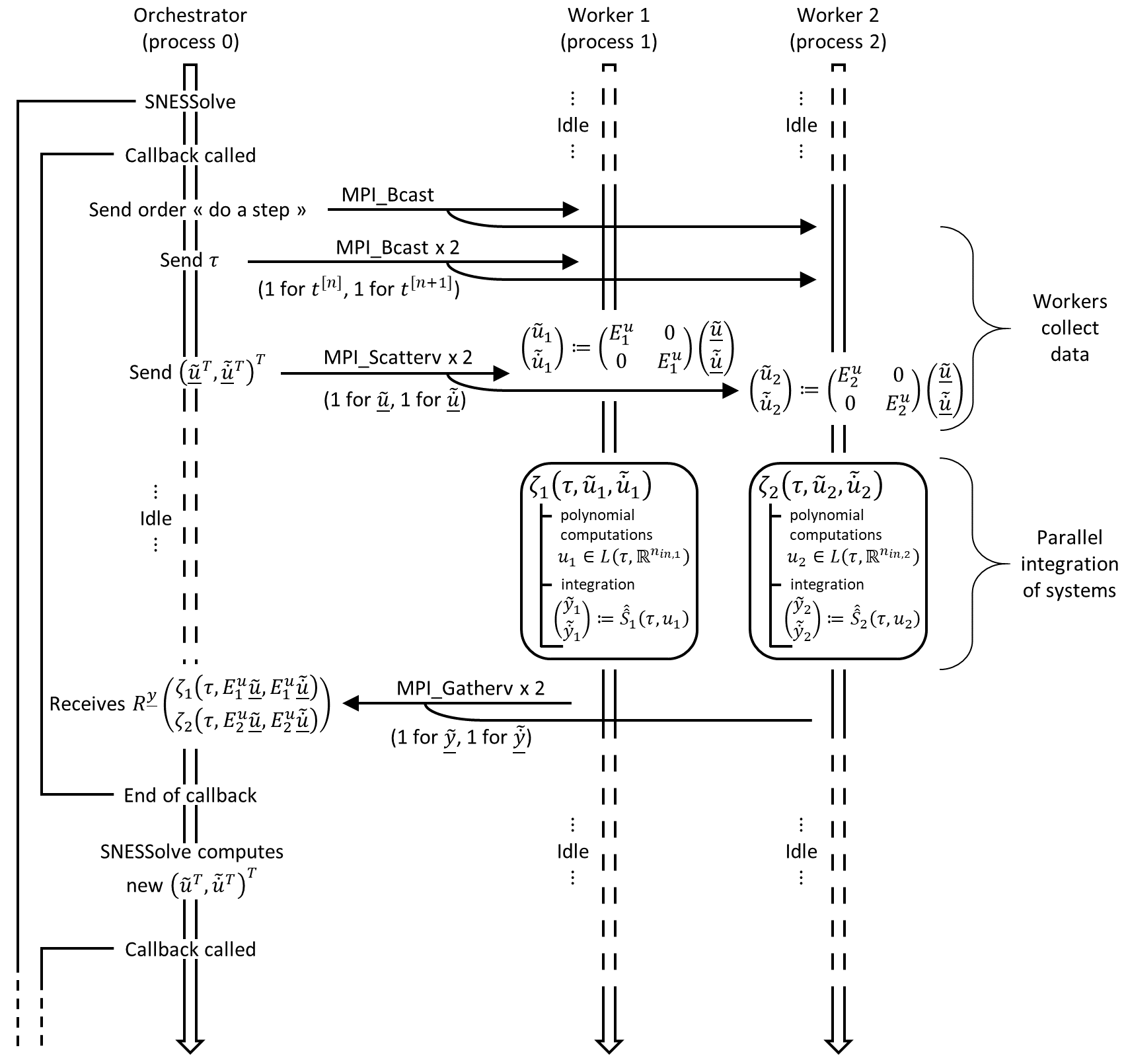}
\captionof{figure}{Workflow of the callback function called by SNESSolve: example with $\nsys=2$ (external first call to the callback is supposed to be already made before \texttt{SNESSolve} is called)}
\label{fig:callback_sheme}
\end{center}

\newpage

\section{Results on a test cases}
\label{section:results_on_a_test_cases}

Difficulties may appear in a co-simulation problem when the coupling is not straightforward. Some of the most difficult cases to solve are the algebraic coupling (addressed in \cite{Gu2004}) arising from causal conflicts, and the multiple feed-through, \ie the case where outputs of a system linearly depend on its inputs, and the connected system(s) have the same behavior. In some case, this may lead to a non-contractant $\Psi_{\tau}$ function. This section presents a test case we designed, belonging to this second category. The fixed-point convergence can be accurately analyzed so that its limitations are highlighted.

\subsection{Test case presentation}
\label{subsection:test_case_presentation}

\begin{center}
\includegraphics[scale=0.35]{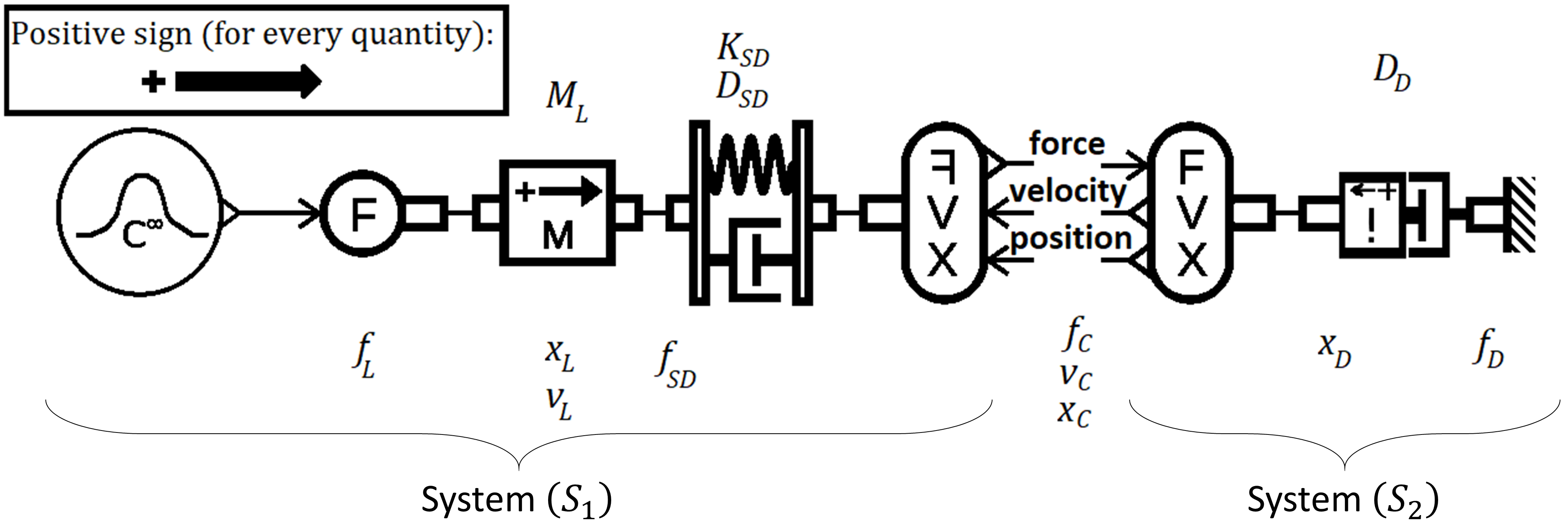}
\captionof{figure}{Mass spring damper with damping reaction modelled with Simcenter Amesim - Parameters are above, variables are below}
\label{fig:test_case_algebraic_loop}
\end{center}

Figure \ref{fig:test_case_algebraic_loop} represents a $1$-mass test case with a classical mechanical coupling on force, velocity and position. These coupling quantities are respectively denoted by $f_c$, $v_c$ and $x_c$. The component on the right represents a damper with a massless plate, computing a velocity (and integrating it to compute a displacement) by reaction to a force input.

We propose the parameters values in table \ref{table:parameters_test_case_algebraic_loop}.

\begin{center}
\captionof{table}{Parameters and initial values of the test case model}
\label{table:parameters_test_case_algebraic_loop}
\vspace{2mm}
\renewcommand{\arraystretch}{1.1}
\begin{tabular}{|l|l|l|}
\hline
Notation & Description & Value \\
\hline
\hline
$M_L$ & Mass of the body in $(S_1)$ & $1$ kg \\
\hline
$K_{SD}$ & Spring rate of the spring in $(S_1)$ & $1$ N/m \\
\hline
$D_{SD}$ & Damper rating of the damper in $(S_1)$ & $1$ N/(m/s) \\
\hline
$D_D$ & Damper rating of the damper in $(S_2)$ & $\in [0.01, 4]$ \\
\hline
\multicolumn{3}{c}{\vspace{-3mm}} \\
\hline
$x_L(0)$ & Initial position of the body in $(S_1)$ & $0$ m \\
\hline
$v_L(0)$ & Initial velocity of the body in $(S_1)$ & $0$ m/s \\
\hline
$x_D(0)$ & Initial position of the plate in $(S_2)$ & $0$ m \\
\hline
\multicolumn{3}{c}{\vspace{-3mm}} \\
\hline
$\tinit$ & Initial time & $0$ s \\
\hline
$\tend$ & Final time & $10$ s \\
\hline
\end{tabular}
\renewcommand{\arraystretch}{1.0}
\end{center}

\noindent All variables will be denoted by either $f$, $v$ or $x$ (corresponding to forces, velocities and positions, respectively) with an index specifying its role in the model (see figure \ref{fig:test_case_algebraic_loop}).

The predefined force $f_L$ is a $C^{\infty}$ function starting from $5$ N and definitely reaching $0$ N at $t=2$ s. The expression of $f_L$ is \eqref{eq:predefined_force_f_L} and the visualization of it is presented on figure \ref{fig:predefined_force_f_L}.

\begin{center}
\includegraphics[scale=0.5]{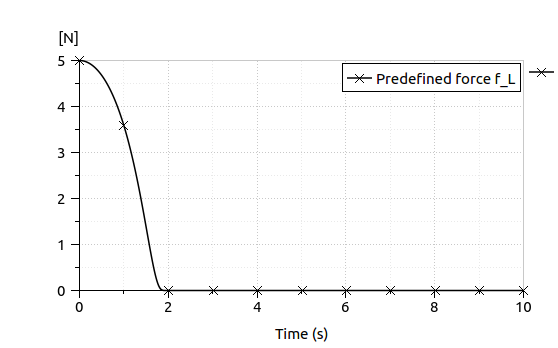}
\captionof{figure}{Predefined force $f_L$}
\label{fig:predefined_force_f_L}
\end{center}

\begin{equation}
\label{eq:predefined_force_f_L}
f_L:\left\{
\begin{array}{lcl}
	[0, 10] & \rightarrow & [0, 5] \\
	t & \mapsto &
		\left\{
		\begin{array}{ll}
		\dspfrac{5}{e^{-1}} e^{\left(\left(\dspfrac{t}{2}\right)^2-\ 1\right)^{-1}} & \text{if}\ t < 2 \\ \\
			0 & \text{if}\ t \geqslant 2 \\
		\end{array}
		\right.
\end{array}
\right.
\end{equation}

\v

\noindent The expected behavior of the model is presented in table \ref{table:main_stages_test_case_model} referring to conventionnal directions of figure \ref{fig:test_case_anim_3D}.

\begin{center}
\includegraphics[scale=0.35]{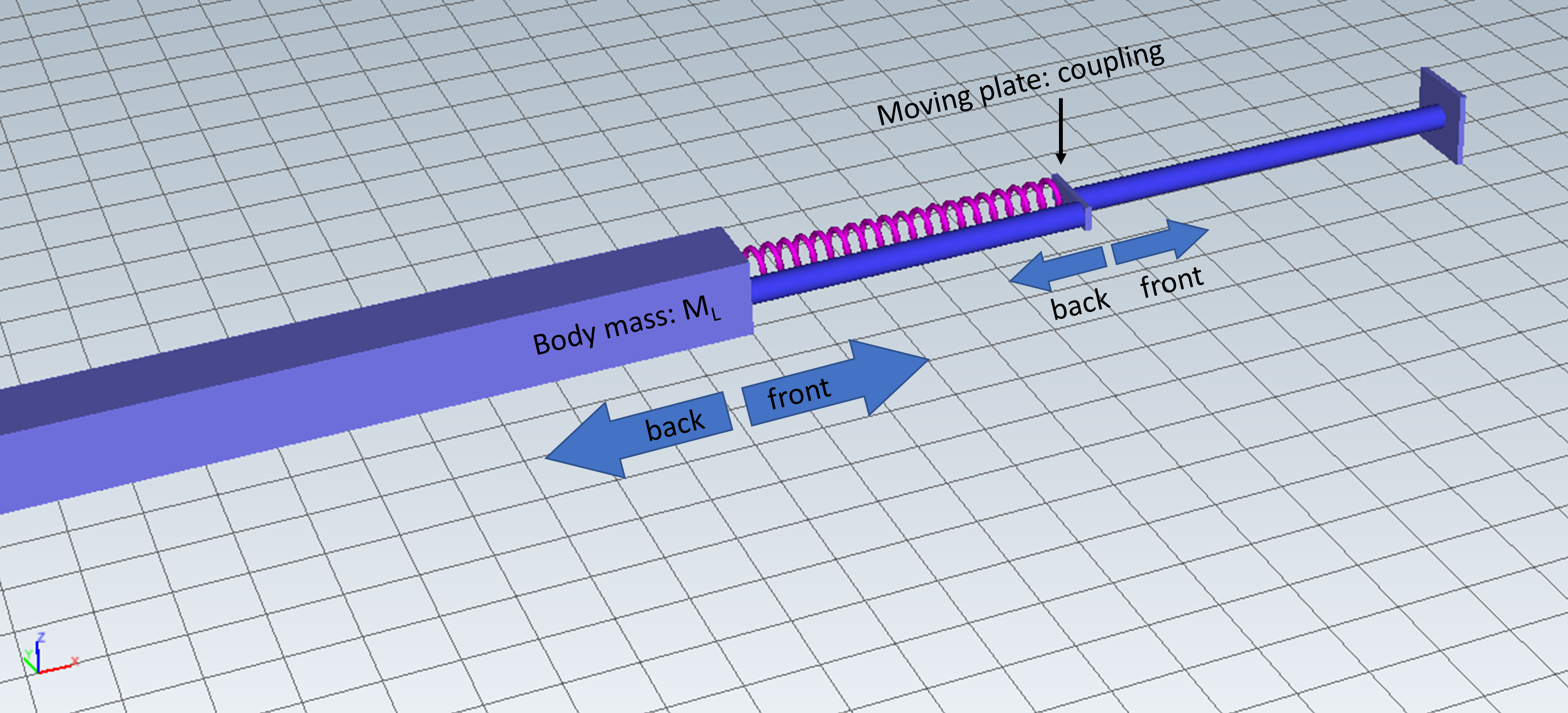}
\captionof{figure}{Test model visualized with Simcenter Amesim}
\label{fig:test_case_anim_3D}
\end{center}

\newpage

\begin{center}
\captionof{table}{Main stages of a simulation of the test case model}
\label{table:main_stages_test_case_model}
\begin{tabular}{|c|c|c|l|}
\hline
\multirow{2}{*}{Stage} & Body & Plate & \multirow{2}{*}{Description} \\
& displacement & displacement & \\
\hline
\hline
$1$ & front & front & Positive $f_L$ pushes everything \\
\hline
$2$ & back & front & \begin{tabular}{l}The spring pushes the body\\backward as it is close to the plate\end{tabular} \\
\hline
$3$ & back & back & \begin{tabular}{l}The spring pulls the plate backwards as\\the body is moving backward with inertia\end{tabular} \\
\hline
$4$ & front & back & \begin{tabular}{l}The spring pulls the body forward as the inertia\\made it go too far in the backward direction\end{tabular} \\
\hline
$5$ & front & front & \begin{tabular}{l}The body is still moving frontward with inertia,\\so the compressed spring pushed the plate forward\end{tabular} \\
\hline
\end{tabular}
\end{center}

\noindent The behavior presented in table \ref{table:main_stages_test_case_model} might slightly change while parameter $D_D$ changes (all other parameters being fixed, see table \ref{table:parameters_test_case_algebraic_loop}).

\subsection{Equations and eigenvalues of the fixed-point callback $\Psi_{\tau}$}
\label{subsection:equations_and_eigenvalues}

Second Newton's law gives:

\begin{equation}
\label{eq:test_case_1}
\begin{array}{lcl}
	\dot{v}_L & = & (f_L + f_{SD}) M_L^{-1} \\
	\dot{x}_L & = & v_L
\end{array}
\end{equation}

\noindent and the spring and damper forces can be expressed the following way:

\begin{equation}
\label{eq:test_case_2}
\begin{array}{lcl}
	f_{SD} & = & K_{SD} ( x_C - x_L ) + D_{SD} (v_C - v_L) \\
	f_C & = & - f_{SD} \\
	f_D & = & - D_D ( 0 - v_C ) \\
	f_C & = & f_D \\
	v_C & = & \nicefrac{f_C}{D_D}
\end{array}
\end{equation}

\noindent leading to the following expressions of the coupled systems:

\begin{equation}
\label{eq:test_case_3}
\begin{array}{l}
	(S_1):\left\{
	\begin{array}{ccl}
		\left(\begin{array}{c}
			\dot{v}_L \\
			\dot{x}_L
		\end{array}\right)
		& = &
		\left(\begin{array}{cc}
			\frac{-D_{SD}}{M_L} & \frac{-K_{SD}}{M_L} \\
			1 & 0
		\end{array}\right)
		\left(\begin{array}{c}
			v_L \\
			x_L
		\end{array}\right)
		+
		\left(\begin{array}{cc}
			\frac{D_{SD}}{M_L} & \frac{K_{SD}}{M_L} \\
			0 & 0
		\end{array}\right)
		\left(\begin{array}{c}
			v_C \\
			x_C
		\end{array}\right)
		+
		\left(\begin{array}{c}
			\frac{f_L}{M_L} \\
			0
		\end{array}\right)
		\\
		f_C & = & \left(D_{SD}\ \ K_{SD}\right)
		\left(\begin{array}{cc}
			v_L \\
			x_L
		\end{array}\right)
		+\left(-D_{SD}\ \ -K_{SD}\right)
		\left(\begin{array}{cc}
			v_C \\
			x_C
		\end{array}\right)
	\end{array}
	\right.
	\\ \\
	(S_2):\left\{
	\begin{array}{ccl}
		\dot{x}_D & = & 0\ x_D + \frac{1}{D_D}\ f_C \\
		\left(\begin{array}{c}
			v_C\\
			x_C
		\end{array}\right)
		& = &
		\left(\begin{array}{c}
			0 \\
			1
		\end{array}\right)
		x_D
		+
		\left(\begin{array}{c}
			\frac{1}{D_D} \\
			0
		\end{array}\right)
		f_C
	\end{array}
	\right.
\end{array}
\end{equation}

\noindent At a given time $t$, we can state the jacobian of $\Psi_{\tau}$ introduced in \eqref{eq:fixed_points_callback} using the expressions of the coupling quantities \eqref{eq:test_case_3}. Indeed, the output variables got at a call are at the same time than the one at which the imposed inputs are reached (end of the macro-step) thanks to the definitions of $\zeta_k$.

\begin{equation}
\label{eq:test_case_jacobian_fixed_point}
J_{\Psi_{\tau}}(\left(\begin{array}{c}f_C\\v_C\\x_c\\\dot{f_C}\\\dot{v_C}\\\dot{x_c}\end{array}\right))
=
\left(\begin{array}{ccccccccc}
	\hhline{-~~~~-~~~}
	\multicolumn{1}{|c|}{0} && -D_{SD} & -K_{SD} && \multicolumn{1}{|c|}{0} && 0 & 0 \\
	\hhline{-~--~-~--}
	\nicefrac{1}{D_{D}} && \multicolumn{1}{|c}{0} & \multicolumn{1}{c|}{0} && 0 && \multicolumn{1}{|c}{0} & \multicolumn{1}{c|}{0} \\
	0 && \multicolumn{1}{|c}{0} & \multicolumn{1}{c|}{0} && 0 && \multicolumn{1}{|c}{0} & \multicolumn{1}{c|}{0} \\
	\hhline{~~--~~~--}
	\vspace{-3mm}
	\\
	\hhline{----~-~~~}
	\multicolumn{4}{|c|}{ } && \multicolumn{1}{|c|}{0} && -D_{SD} & -K_{SD} \\
	\hhline{~~~~~-~--}
	\multicolumn{4}{|c|}{\text{Block}} && \nicefrac{1}{D_{D}} && \multicolumn{1}{|c}{0} & \multicolumn{1}{c|}{0} \\
	\multicolumn{4}{|c|}{ } && 0 && \multicolumn{1}{|c}{0} & \multicolumn{1}{c|}{0} \\
	\hhline{----~~~--}
\end{array}\right)
\renewcommand{\arraystretch}{1.0}
\end{equation}

\noindent The framed zeros are "by-design" zeros: indeed, systems never produce outputs depending on inputs given to other systems. The block called "Block" in \eqref{eq:test_case_jacobian_fixed_point} depends on the methods used to retrieve the time-derivatives of the coupling quantities (see \eqref{eq:step_function_pratical_extended} and its finite differences version). Nevertheless, this block does not change the eigenvalues of $J_{\Psi_{\tau}}$ as it is a block-triangular matrix. Indeed, the characteristic polynomial of $I_6 - \lambda J_{\Psi_{\tau}}$ is the product of the determinant of the two $3\times 3$ blocks on the diagonal of $I_6 - \lambda J_{\Psi_{\tau}}$. The eigenvalues of $J_{\Psi}$ are:

\begin{equation}
\label{eq:test_case_eigenvalues}
0,\ +1\!\textit{i}\sqrt{\frac{D_{SD}}{D_D}},\ -1\!\textit{i}\sqrt{\frac{D_{SD}}{D_D}}\hspace{5mm}\text{(each with a multiplicity of $2$)}
\end{equation}

\noindent Hence, the following relation between the parameters and the spectral radius can be shown (given $D_D > 0$ and $D_{SD} = 1 > 0$):

\begin{equation}
\label{eq:test_case_spectral_radius}
\varrho\left(J_{\Psi_{\tau}}\right)\left\{
\begin{array}{lcl}
	< 1 & \ \text{if}\ D_{SD} < D_D \\
	\geqslant 1 & \ \text{if}\ D_{SD} \geqslant D_D
\end{array}
\right.
\end{equation}

\noindent We can thus expect that the classical IFOSMONDI co-simulation algorithm based on a fixed-point method \cite{Eguillon2019} cannot converge on this model when the damping ratio of the component on the right of the model (see figure \ref{fig:test_case_algebraic_loop}) is smaller than the damping ratio of the spring-damper component.

We will process several simulations with different values of $D_D$ leading to different values of $\varrho(J_{\Psi_{\tau}})$. These values and the expected movement of the body of the system is plotted in figure \ref{fig:solutions_for_various_rho}.

\begin{center}
\includegraphics[scale=0.35]{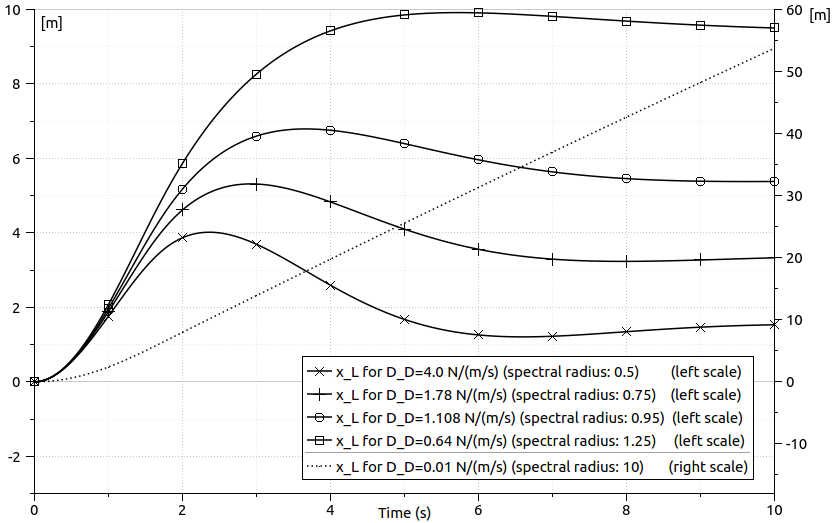}
\captionof{figure}{Displacement of the mass ($x_L$) for different damping ratios of the right damper ($D_D$) simulated on a monolithic model (without co-simulation). Associated spectral radii of $J_{\Psi}$ are recalled for futher coupled formulations.}
\label{fig:solutions_for_various_rho}
\end{center}

\subsection{Results}
\label{subsection:results}

As the PETSc library enables to easily change the parameters of the JFM (as explained in subsection \ref{subsection:using_petsc_for_the_jfm}), three methods have been used in the simulations:

\begin{itemize}
\ite NewtonLS: a Newton based non-linear solver that uses a line search,
\ite Ngmres: the non-linear generalized minimum residual method \cite{Oosterlee2000}, and
\ite Anderson: the Anderson mixing method \cite{Anderson1965}
\end{itemize}

First of all, simulations have been processed with all these JFMs (with parameters exhaustively defined in appendix \ref{appendix:parameters_of_the_petsc_non_linear_solvers}) within IFOSMONDI-JFM, the classical IFOSMONDI algorithm (denoted hereafter as "Fixed-point"), and the original explicit zero-order hold co-simulation method (sometimes referred to as non-iterative Jacobi). The error is defined as the mean of the normalized $L^2$ errors on each state variable of both systems on the whole $[\tinit, \tend]$ domain. The reference is the monolithic simulation (of the non-coupled model) done with Simcenter Amesim. Such errors are presented for a contractant case ($D_D=4$ N, so $\varrho(J_{\Psi_{\tau}})=0.5$) in figure \ref{fig:err_vs_dt_contractant}. For a non-contractant case ($D_D=0.64$ N, so $\varrho(J_{\Psi_{\tau}})=1.25$), analog plots are presented in figure \ref{fig:err_vs_dt_noncontractant}.

\begin{center}
\includegraphics[scale=0.43]{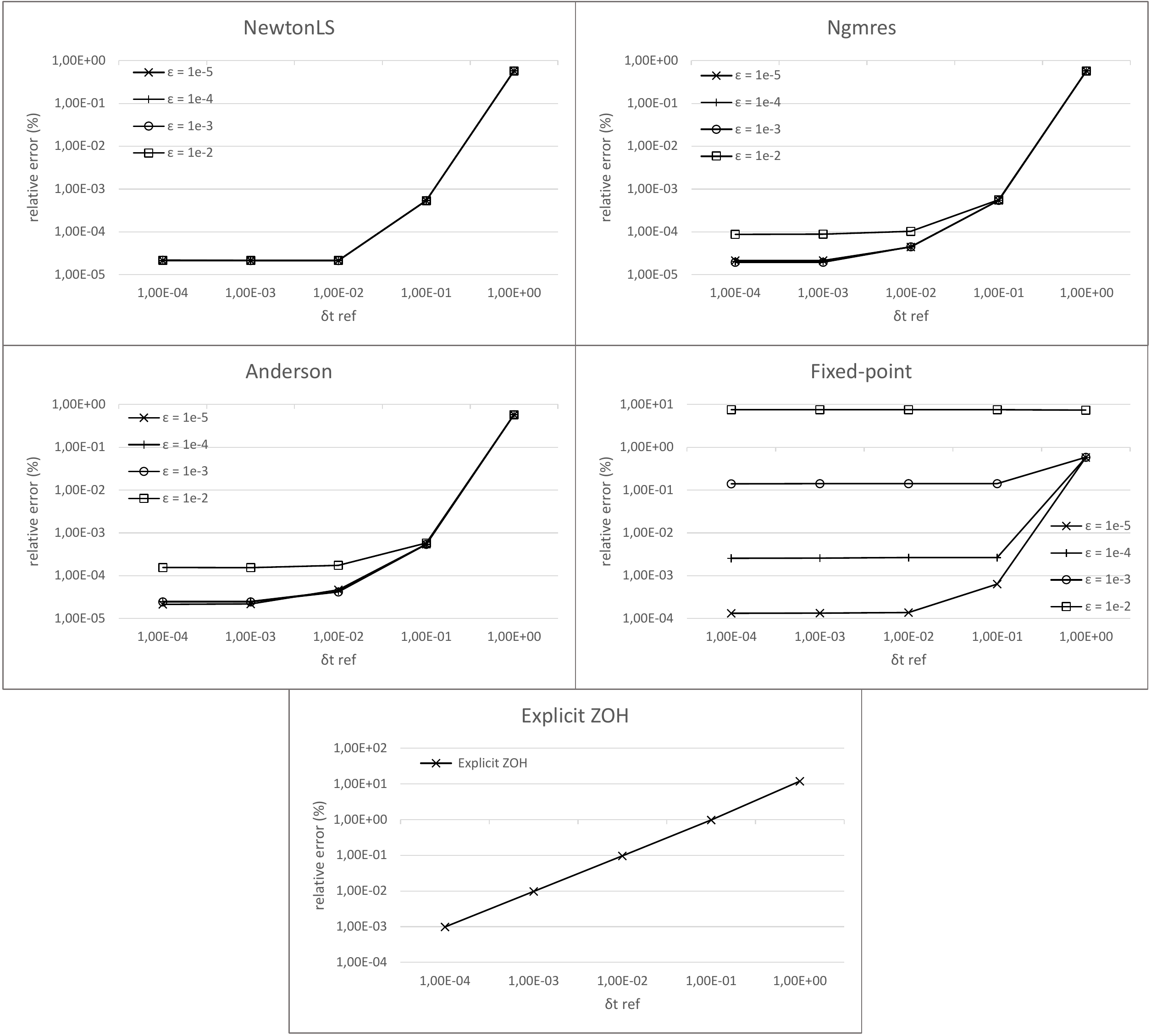}
\vspace{-3mm}
\captionof{figure}{Error accross $\dtref$ with different methods on a contractant case ($D_D=4.0$, $\rho(J_{Psi})=0.5$) - NewtonLS, Ngmres and Anderson are matrix-free iterative methods used with the IFOSMONDI-JFM algorithm, Fixed-point is the classical IFOSMONDI algorithm, and Explicit ZOH is the non-iterative zero-order hold fixed-step co-simulation}
\label{fig:err_vs_dt_contractant}
\end{center}

\begin{center}
\includegraphics[scale=0.45]{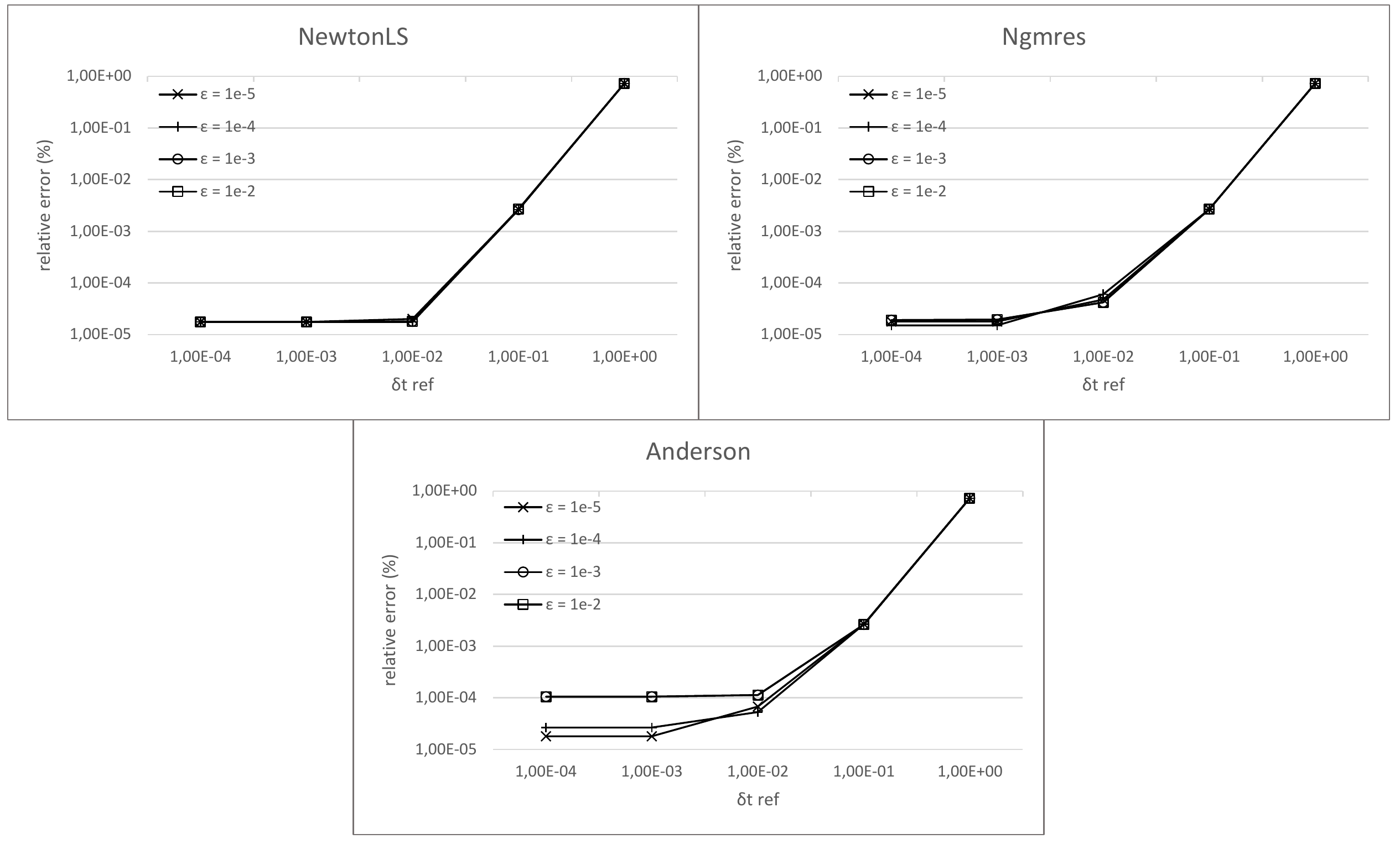}
\captionof{figure}{Error accross $\dtref$ with different methods on a non-contractant case ($D_D=0.64$, $\rho(J_{Psi})=1.25$) - NewtonLS, Ngmres and Anderson are matrix-free iterative methods used with the IFOSMONDI-JFM algorithm}
\label{fig:err_vs_dt_noncontractant}
\end{center}

\noindent As expected, the simulations failed (diverged) with fixed-point method (classical IFOSMONDI) for the non-contractant case. Moreover, the values given to the system were too far from physically-possible values with the explicit ZOH co-simulation algorithm, so the internal solvers of systems $(S_1)$ and $(S_2)$ failed to integrate. These are the reason why these two methods lead to no curve on figure \ref{fig:err_vs_dt_noncontractant}.

Nonetheless, the three versions of IFOSMONDI-JFM algorithm keep producing reliable results with an acceptable relative error (less than $1\%$) when $\dtref \geqslant 0.1$ s.
\\

On figures \ref{fig:err_vs_dt_contractant} and \ref{fig:err_vs_dt_noncontractant}, IFOSMONDI-JFM method seems to solve the problem with a good accuracy regardless of the value of the damping ratio $D_D$. In order to confirm that, several other values have been tried: the ones for which the solution has been computed and plotted in figure \ref{fig:solutions_for_various_rho}. The error is presented, but also the number of iterations and the number of integrations (calls to $\zeta_k$, i.e. calls to $\gamma_{\tau}$ for IFOSMONDI-JFM or to $\Psi_{\tau}$ for classical IFOSMONDI). Although for the fixed-point case (classical IFOSMONDI) the number of iteration is the same than the number of integration, for the IFOSMONDI-JFM algorithm the number of iterations is the one of the underlying non-linear solver (NewtonLS, Ngmres or Anderson), and there might be a lot more integrations than iterations of the non-linear method. These results are presented in figure \ref{fig:err_nbIter_nbIntegr_across_rho}.

\begin{center}
\includegraphics[scale=0.45]{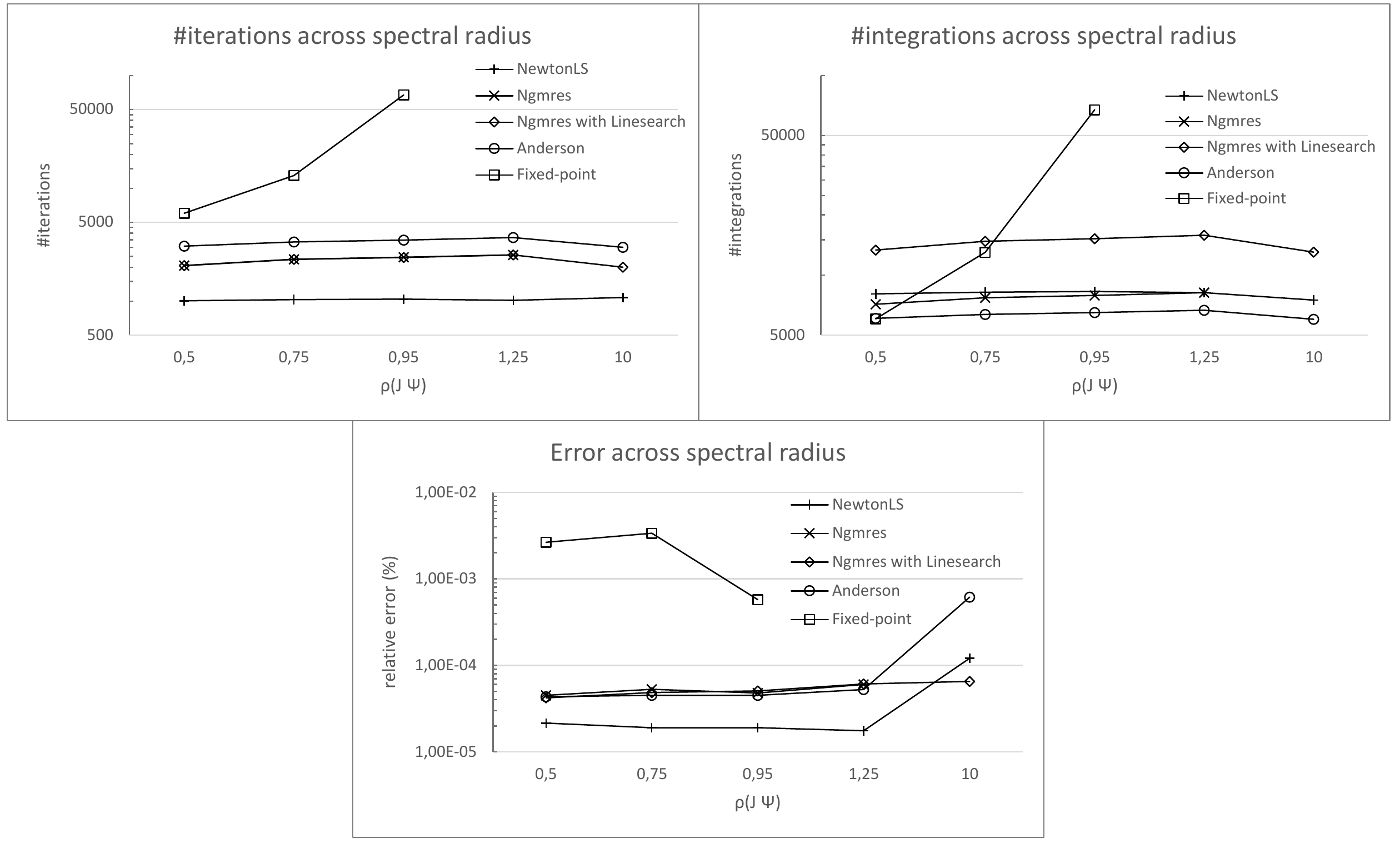}
\captionof{figure}{Total number of iterations, integrations, and error across spectral radius of $J_{\Psi}$ for different methods (Fixed-point corresponds to classical IFOSMONDI algorithms, and all other methods are used with the IFOSMONDI-JFM version). All co-simulation ran with $\varepsilon = 10^{-4}$ and $\dtref = 10^{-2}$}
\label{fig:err_nbIter_nbIntegr_across_rho}
\end{center}

\noindent As expected, the threshold of $\varrho(J_{\Psi_{\tau}})=1$ (\ie $D_D=D_{SD}=1$) is critical for the fixed-point method. The IFOSMONDI-JFM method not only can overpass this threshold, but no significant extra dificulty appears to solve the problem in the non-contractant cases, except for the Ngmres non-linear solver (which failed to converge with $D_D=0.01$, so with $\varrho(J_{\Psi_{\tau}})=10$). However, regarding the Ngmres method, the variant that uses line search converges in all cases. Eventhough the latter requires more integrations than other JFMs, it is more robust to high values of $\varrho(J_{\Psi_{\tau}})$. The parameters of this line search are detailed on table \ref{table:parameters_of_the_ngmres_with_linesearch_method} in appendix \ref{appendix:parameters_of_the_petsc_non_linear_solvers}.

The NewtonLS and Anderson methods show a slightly bigger error on this "extreme" case of $\varrho(J_{\Psi_{\tau}})=10$, yet it stays under $0.001\%$ which is completely acceptable.

Among those two JFMs (NewtonLS and Anderson), the trend that can be observed on figure \ref{fig:err_nbIter_nbIntegr_across_rho} shows that NewtonLS is always more accurate than Anderson, yet it always requires a bigger amount of integrations. We can stand that IFOSMONDI-JFM is more \textit{accuracy-oriented} on this model when it is based on the NewtonLS JFM, and more \textit{speed-oriented} on this model when it is based on the Anderson JFM (for the same $\dtref$ and $\varepsilon$). For high values of $\varrho(J_{\Psi_{\tau}})$, \textit{accuracy-oriented} simulations are achieved thanks to the Ngmres JFM with line search more than the NewtonLS one.

Finally, smaller errors are obtained with IFOSMONDI-JFM and with less iterations than classical IFOSMONDI. Yet, the time consumption is directly linked with the number of integrations, not with the number of iterations of the underlying non-linear solver. The total number of integrations does not increase across the problem difficulty (increasing with $\varrho(J_{\Psi_{\tau}})$), and the non-linear methods within IFOSMONDI-JFM do not even require more integrations that the fixed-point one for most of the values of $D_D$ for which the classical IFOSMONDI algorithm does not fail.

\section{Conclusion}
\label{section:conclusion}

IFOSMONDI-JFM method not only enables to solve problems that were impossible to solve with the classical IFOSMONDI method, it also requires less iterations to converge on the test case of section \ref{section:results_on_a_test_cases} when the parameterization enables both methods to solve the problem.

Despite a number of integration greater than one for every iteration (contrary to the classical IFOSMONDI algorithm), IFOSMONDI-JFM does not require a lot more integrations than classical IFOSMONDI. In most of the cases, for the considered test case, IFOSMONDI-JFM even requires less integrations than classical IFOSMONDI, and the resulting solution is always more accurate (for the same $\dtref$ and $\varepsilon$). The matrix-free aspect of the underlying solvers used with IFOSMONDI-JFM are one of the causes of the small amount of integrations.

The IFOSMONDI-JFM algorithm takes advantages from the $C^1$ smoothness of classical IFOSMONDI algorithm \cite{Eguillon2019} without the delay it implies in \cite{Busch2019} (thanks to its iterative aspect), the coupling constraint is satisfied both at left and right of every communication time thanks to the underlying non-linear solvers of PETSc \cite{PetscWebPage}. The iterative part does not need a finite differences estimation of the jacobian matrix like in \cite{Schweizer2015} or a reconstruction of it like in \cite{Sicklinger2014}.

The resulting algorithm even solves co-simulation problems for which the fixed-point formulation would involve a non-contractant coupling function $\Psi_{\tau}$.

Finally, the test case introduced in \ref{subsection:test_case_presentation} can be reused to test the robustness of various co-simulation methods as the model is relatively simple and the difficulty can easily be increased or decreased in a quantifiable way.

\renewcommand{\lstlistingname}{Listing}

\bibliographystyle{splncs04}
\bibliography{AES_refs}

\newpage

%
%


\begin{subappendices}
\renewcommand{\thesection}{\Alph{section}}

\section{Parameters of the PETSc non-linear solvers}
\label{appendix:parameters_of_the_petsc_non_linear_solvers}

The JFMs mentionned in this document (see definition in \ref{subsection:a_word_on_jfm_accronym}) refer to PETSc non-linear solvers, so-called '\texttt{SNES}' in the PETSc framework.

The parameters of these methods where the default one, except the explicitely mentionned ones. The following tables recaps these options. For furthe definition of their meaning, see \cite{PetscWebPage,Anderson1965,Oosterlee2000}.

\begin{center}
\captionof{table}{Parameters of the NewtonLS method}
\label{table:parameters_of_the_newtonls_method}
\renewcommand{\arraystretch}{1.5}
\begin{tabular}{|c|m{7cm}|c|}
\hline
\multicolumn{1}{|>{\centering}m{4cm}|}{PETSc argument: \begin{footnotesize}\texttt{-snes\_linesearch\_<...>}\end{footnotesize}} & \multicolumn{1}{c|}{Description} & Value \\
\hline
\hline
\verb!type! & Select line search type & \verb!bt! \\
\hline
\verb!order! & Selects the order of the line search for bt & $3$ \\
\hline
\verb!norms! & Turns on/off computation of the norms for basic line search & \verb!TRUE! \\
\hline
\verb!alpha! & Sets alpha used in determining if reduction in function norm is sufficient & $0.0001$ \\
\hline
\verb!maxstep! & Sets the maximum stepsize the line search will use & $10^8$ \\
\hline
\verb!minlambda! & Sets the minimum lambda the line search will tolerate & $10^{-12}$ \\
\hline
\verb!damping! & Damping factor used for basic line search & $1$ \\
\hline
\verb!rtol! & Relative tolerance for iterative line search & $10^{-8}$ \\
\hline
\verb!atol! & Absolute tolerance for iterative line search & $10^{-15}$ \\
\hline
\verb!ltol! & Change in lambda tolerance for iterative line search & $10^{-8}$ \\
\hline
\verb!max_it! & Maximum iterations for iterative line searches & $40$ \\
\hline
\verb!keeplambda! & Use previous lambda as damping & \verb!FALSE! \\
\hline
\verb!precheck_picard! & Use a correction that sometimes improves convergence of Picard iteration & \verb!FALSE! \\
\hline
\end{tabular}
\renewcommand{\arraystretch}{1.0}
\end{center}

\newpage

\begin{center}
\captionof{table}{Parameters of the Anderson method}
\label{table:parameters_of_the_anderson_method}
\renewcommand{\arraystretch}{1.5}
\begin{tabular}{|c|m{7cm}|c|}
\hline
\multicolumn{1}{|>{\centering}m{4cm}|}{PETSc argument: \begin{footnotesize}\texttt{-snes\_anderson\_<...>}\end{footnotesize}} & \multicolumn{1}{c|}{Description} & Value \\
\hline
\hline
\verb!m! & Number of stored previous solutions and residuals & $30$ \\
\hline
\verb!beta! & Anderson mixing parameter & $1$ \\
\hline
\verb!restart_type! & Type of restart & \verb!NONE! \\
\hline
\verb!restart_it! & Number of iterations of restart conditions before restart & $2$ \\
\hline
\verb!restart! & Number of iterations before periodic restart & $30$ \\
\hline
\end{tabular}
\renewcommand{\arraystretch}{1.0}
\end{center}

\begin{center}
\captionof{table}{Parameters of the Ngmres method (not Ngmres with line search)}
\label{table:parameters_of_the_ngmres_method}
\renewcommand{\arraystretch}{1.5}
\begin{tabular}{|c|m{7cm}|c|}
\hline
\multicolumn{1}{|>{\centering}m{4cm}|}{PETSc argument: \begin{footnotesize}\texttt{-snes\_ngmres\_<...>}\end{footnotesize}} & \multicolumn{1}{c|}{Description} & Value \\
\hline
\hline
\verb!select_type! & Choose the select between candidate and combined solution & \verb!DIFFERENCE! \\
\hline
\verb!restart_type! & Choose the restart conditions & \verb!DIFFERENCE! \\
\hline
\verb!candidate! & Use NGMRES variant which combines candidate solutions instead of actual solutions & \verb!FALSE! \\
\hline
\verb!approxfunc! & Linearly approximate the function & \verb!FALSE! \\
\hline
\verb!m! & Number of stored previous solutions and residuals & $30$ \\
\hline
\verb!restart_it! & Number of iterations the restart conditions hold before restart & $2$ \\
\hline
\verb!gammaA! & Residual tolerance for solution select between the candidate and combination & $2$ \\
\hline
\verb!gammaC! & Residual tolerance for restart & $2$ \\
\hline
\verb!epsilonB! & Difference tolerance between subsequent solutions triggering restart & $0.1$ \\
\hline
\verb!deltaB! & Difference tolerance between residuals triggering restart & $0.9$ \\
\hline
\verb!single_reduction! & Aggregate reductions & \verb!FALSE! \\
\hline
\verb!restart_fm_rise! & Restart on residual rise from x\_M step & \verb!FALSE! \\
\hline
\end{tabular}
\renewcommand{\arraystretch}{1.0}
\end{center}

\begin{center}
\captionof{table}{Parameters of the Ngmres with linsearch method}
\label{table:parameters_of_the_ngmres_with_linesearch_method}
\renewcommand{\arraystretch}{1.5}
\begin{tabular}{|c|m{7cm}|c|}
\hline
\multicolumn{1}{|>{\centering}m{4cm}|}{PETSc argument: \begin{footnotesize}\texttt{-snes\_ngmres\_<...>}\end{footnotesize}} & \multicolumn{1}{c|}{Description} & Value \\
\hline
\hline
\verb!select_type! & Choose the select between candidate and combined solution & \verb!LINESEARCH! \\
\hline
\multicolumn{3}{c}{\vdots} \\
\multicolumn{3}{c}{All other options of table \ref{table:parameters_of_the_ngmres_method} are the same} \\
\multicolumn{3}{c}{\vdots} \\
\hline
\multicolumn{1}{|>{\centering}m{4cm}|}{PETSc argument: \begin{footnotesize}\texttt{-snes\_linesearch\_<...>}\end{footnotesize}} & \multicolumn{1}{c|}{Description} & Value \\
\hline
\hline
\verb!type! & Select line search type & \verb!basic! \\
\hline
\verb!order! & Selects the order of the line search for bt & $0$ \\
\hline
\verb!norms! & Turns on/off computation of the norms for basic linesearch & \verb!TRUE! \\
\hline
\verb!maxstep! & Sets the maximum stepsize the line search will use & $10^8$ \\
\hline
\verb!minlambda! & Sets the minimum lambda the line search will tolerate & $10^{-12}$ \\
\hline
\verb!damping! & Damping factor used for basic line search & $1$ \\
\hline
\verb!rtol! & Relative tolerance for iterative line search & $10^{-8}$ \\
\hline
\verb!atol! & Absolute tolerance for iterative line search & $10^{-15}$ \\
\hline
\verb!ltol! & Change in lambda tolerance for iterative line search & $10^{-8}$ \\
\hline
\verb!max_it! & Maximum iterations for iterative line searches & $1$ \\
\hline
\verb!keeplambda! & Use previous lambda as damping & \verb!FALSE! \\
\hline
\verb!precheck_picard! & Use a correction that sometimes improves convergence of Picard iteration & \verb!FALSE! \\
\hline
\end{tabular}
\renewcommand{\arraystretch}{1.0}
\end{center}

\end{subappendices}

\end{document}